\def\hybrid{\topmargin 0pt      \oddsidemargin 0pt
        \headheight 0pt \headsep 0pt
        \textwidth 160true mm       
        \textheight 231true mm         
        \marginparwidth 0.0in
        \parskip 0pt plus 1pt   \jot = 1.5ex}
\catcode`\@=11
\def\marginnote#1{}

\newcount\hour
\newcount\minute
\newtoks\amorpm
\hour=\time\divide\hour by60
\minute=\time{\multiply\hour by60 \global\advance\minute by-\hour}
\edef\standardtime{{\ifnum\hour<12 \global\amorpm={am}%
        \else\global\amorpm={pm}\advance\hour by-12 \fi
        \ifnum\hour=0 \hour=12 \fi
        \number\hour:\ifnum\minute<10 0\fi\number\minute\the\amorpm}}
\edef\militarytime{\number\hour:\ifnum\minute<10 0\fi\number\minute}

\def\draftlabel#1{{\@bsphack\if@filesw {\let\thepage\relax
   \xdef\@gtempa{\write\@auxout{\string
      \newlabel{#1}{{\@currentlabel}{\thepage}}}}}\@gtempa
   \if@nobreak \ifvmode\nobreak\fi\fi\fi\@esphack}
        \gdef\@eqnlabel{#1}}
\def\@eqnlabel{}
\def\@vacuum{}
\def\draftmarginnote#1{\marginpar{\raggedright\scriptsize\tt#1}}

\def\draft{\oddsidemargin -.5truein
        \def\@oddfoot{\sl preliminary draft \hfil
        \rm\thepage\hfil\sl\today\quad\militarytime}
        \let\@evenfoot\@oddfoot \overfullrule 3pt
        \let\label=\draftlabel
        \let\marginnote=\draftmarginnote
   \def\@eqnnum{(\theequation)\rlap{\kern\marginparsep\tt\@eqnlabel}%
\global\let\@eqnlabel\@vacuum}  }

\relax
\hybrid

\documentclass[10pt]{article}
\usepackage{amsmath}
\usepackage{amsfonts}
\usepackage{amssymb,amscd}

\def\ba{\begin{eqnarray}}
\def\ea{\end{eqnarray}}

\def\lb{\label}
\def\Tr{{\rm Tr}}
\def\be{\begin{equation}}
\def\ee{\end{equation}}
\def\qed{\rule{5pt}{5pt}}

\newcommand{\bea}{\begin{eqnarray}}
\newcommand{\eea}{\end{eqnarray}}

\def\theequation{\thesection.\arabic{equation}}

\begin{document}

{}

\vspace{2cm}

\begin{center}
{\Large \bf Bethe subalgebras in affine Birman--Murakami--Wenzl algebras and
flat connections  for $q$-KZ equations.}6
\end{center}

\vspace{0.5cm}

\begin{center}
{\bf A.P.Isaev$^{*}$, A.N.Kirillov$^{**}$ and V.O.Tarasov$^{***}$}
\end{center}

\vspace{1cm}

\begin{center}
$^*$ Bogoliubov Laboratory of Theoretical Physics, JINR, \\
 141980, Dubna, Moscow region, \\
 and ITPM, M.V.Lomonosov Moscow State University, Russia \\
 E-mail: isaevap@theor.jinr.ru
\end{center}

\begin{center}
$^{**}$ Research Institute of Mathematical Sciences, RIMS, \\
Kyoto University, Sakyo-ku, 606-8502, Japan \\
{\it URL: ~~~~http://www.kurims.kyoto-u.ac.jp/\kern-.05cm$\tilde{\quad}$\kern-.17cm kirillov } \\
and \\
The Kavli Institute for the Physics and Mathematics of the Universe
( IPMU ),\\ 5-1-5 Kashiwanoha,  Kashiwa, 277-8583, Japan \\
and \\
Department of Mathematics, National Research University Higher School of
Economics, \\
117312, Moscow, Vavilova str. 7,~ Russia
\end{center}

\begin{center}
$^{***}$Department of Mathematical Sciences, \\
Indiana University -- Purdue University Indianapolis \\
402 North Blackford St, Indianapolis, IN 46202-3216, USA \\
and  \\
St. Petersburg Branch of Steklov Mathematical Institute
Fontanka 27, \\ St. Petersburg, 191023,~ Russia

\vspace{5mm}

\begin{flushright}
\begin{minipage}{120mm}
\it Dedicated to Professor Rodney Baxter  on the occasion of his 75th Birthday.
\end{minipage}
\end{flushright}

\vspace{6mm}

\end{center}

\vspace{1cm}

{\bf Abstract.} Commutative sets of Jucys--Murphy
elements for affine braid groups of $A^{(1)},B^{(1)},C^{(1)},D^{(1)}$
types were defined.
Construction of $R$-matrix representations of the affine braid
group of type $C^{(1)}$ and  its distinguish commutative subgroup generated
by the $C^{(1)}$-type  Jucys--Murphy elements are given. We describe a general method to produce flat connections for the two-boundary quantum Knizhnik-
Zamolodchikov equations as necessary  conditions for Sklyanin's type
transfer matrix associated with the two-boundary multicomponent Zamolodchikov
algebra to be invariant under the action of the $C^{(1)}$-type Jucys--Murphy
elements. We specify our general construction to the case of the Birman--Murakami--Wenzl algebras ($BMW$ algebras for short). As an application we suggest a
baxterization of the Dunkl--Cherednik elements $Y'$s in the double affine
Hecke algebra of type $A$.

\vspace{1cm}


\vspace{0.5cm}


\vspace{1cm}
\section{Introduction}
\setcounter{equation}0

The quantum Knizhnik--Zamolodchikov equation (q-KZ equation for shot) is a
system of difference equations which   has been introduced by F.Smirnov
\cite{Sm}, \cite{Sm2},  during the study of form factors of integrable models,
 and independently, by  I. Frenkel and N.Reshetikhin,  \cite{FreRe} during the
study of the representation theory of quantum affine algebras. Since that
time the literature that enter into the treatment of qKZ equations, their
generalizations and  applications, are enormous. We mention here only a few:

$\bullet$~~\cite{JMN},~ which is concerned to the study of correlation
functions of  integrable systems;

$\bullet$~~\cite{Ch},~which is devoted to applications to the representation
theory of affine Hecke algebras;

$\bullet$~~\cite{KZ}, \cite{Z}, \cite{PZ},~~which are concerned to the study of variety
 applications to Algebraic Combinatorics and Algebraic Geometry of certain
class of solutions to (boundary) q-KZ equations.

$\bullet$~~\cite{RSV}, ~devoted to the study of Jackson integral solutions of
the boundary quantum Knizhnik-Zamolodchikov equation(s) with applications to
the representation theory of quantum affine algebra  $U_{q}(\widehat{{\mathfrak{sl}}(2)})$.

In the present paper we describe a general method for construction of
{\it two-boundary} quantum KZ equations associated with affine Birman--Murakami--Wenzel algebras (BMW algebras) \cite{BirWen}, \cite{Mur00}, \cite{We},
\cite{DRV},~  and give several examples to illustrate our method. The
underlying idea of our
construction is to describe relations/equations among the generators of the
multicomponent
 two-boundary Zamolodchikov algebras \cite{GoshZam} which \underline{imply} that the natural
action of the distinguish commutative subgroup of the affine braid group
$B_n(C^{(1)})$ of type $C^{(1)}$ generated by the Jucys--Murphy elements
$\{ JM_{i} \},~i=1,2,\ldots, n,$~~
~\underline{preserves} the ~``monodromy matrix''~ associated with the
Zamolodchikov algebras in question, see Sections 2, 3 and 4  for details. For
example, in Section 2 we describe {\it distinguish} commutative subgroups in
the (non-twisted) affine braid groups of classical types. The
generators of these distinguish subgroups will be called {\it universal Jucys--
Murphy elements}, or  $JM$-elements for short. Note that the well-known
$JM$-elements in the group ring of the symmetric group \cite{Ju0}, or (affine) Hecke  \cite{GiNi}, \cite{Doi},
Birman--Murakami-Wenzl \cite{IsOg} and  cyclotomic Hecke (and cyclotomic BMW) algebras,
are images of the universal $JM$-elements. The main objective of our paper is
to construct {\it Baxterization} of
the $JM$-elements in the affine  Birman--Murakami--Wenzl  algebras of type
$C^{(1)}$, i.e. to construct mutually commuting family of elements $JM_i(x) \in BMW(C^{1)}) \otimes \mathbb{Q}(x)$ depending on spectral parameter $x$, such
that $JM_i(0)=JM_i,~ \forall i$.  \\

Now let us say few words about the content of our paper. \\
As it was mentioned, in \underline{Section~2} we recall definitions of
{\it distinguish} commutative subgroups in the affine braid groups of classical types. Since the generators of these commutative subgroups are the major origin of the Jucys--Murphy elements in a big variety of algebras, we include the
definitions and  proofs of universal $JM$-elements basic properties.

We want to \underline{stress} that in all known cases, such as the group ring
of the symmetric groups, (affine, cyclotomic) Hecke, Brauer, $BMW$ algebras, the corresponding $JM$-elements come from the distinguish commutative subgroup in
the corresponding (affine) braid group of classical type. In fact, birational
representations of affine braid group associated with semisimple Lie algebras,  give rise to the well-known and widely used integrable systems such as
Heisenberg chains and Gaudin models, \cite{IsOg2}, \cite{IK},
 Painlev\'{e} equations, \cite{No} and the literature quoted therein.   \\

In \underline{Section~\ref{sec3}} we describe a way how to construct $R$-matrix
representations of the affine braid group $B_n(C^{(1)})$ of type $C^{(1)}$, and
 use these constructions to define the corresponding {\it quantum $qKZ$ equations}  and  two sets of {\it flat connections} associated with the former. \\

\underline{Section \ref{sec4}} contains one of our main results  concerning of construction of {\it flat connections} based on the study of two-boundary (multi-component) Zamolodchikov algebras. Namely, $qKZ$ equations are making their appearance to  ensure that the two boundary Zamolodchikov algebra in question is invariant under the action of the distinguish commutative subgroup in the corresponding affine braid group. In \underline{Section \ref{sec4}.2} we present our main construction, namely that of {\it flat connections} for quantum Knizhnik--Zamolodchikov equations derived from the study of two boundary Zamolodchikov algebra and the
$B_n(C^{(1)})$ universal Jucys--Murphy elements.

In \underline{Section \ref{aBMW}} we specify our general constructions presented in Section \ref{sec4} to the case of affine  $BMW$ algebras, and construct flat connections for
the algebra $BMW(C^{(1)})$. To pass from general construction to the case of
the affine Birman--Murakami--Wenzl algebras of type $C^{(1)}$, we rely on the
use of embedding the braid group $B_n(C^{(1)})$ into the algebra $BMW(C^{(1)})$.

In \underline{Section \ref{sec6}} we construct {\it baxterized Jucys--Murphy elements}
in the affine $BMW$ algebras. Our approach is based on Sklyanin's transfer
 matrix method
\footnote{ Naive replacement of generators $T_i$  in $(\ref{jucys1})$  by its baxterization $T_i(u/v)$ defined in (\ref{a00}), leads to the set of elements in the $BMW$
algebra,  which do not commute in general}, \cite{Skl},\cite{Skly}.
 The key to apply  the Sklyanin transfer matrix method to construction of
{\it baxterized  $JM$-elements $\bar{y}_n(x;\vec{z}_{(n)})$}, see (\ref{xxz55}),
~\underline{lies} in the fact that the family of algebras $\{ BMW_n(C) \}_{n \ge 1}$ can be provided with the {\it Markov trace}, namely, there exists a
unique homomorphism
$$ Tr_{n+1} :  BMW_{n+1}(C) \longrightarrow BMW_n(C),~~\forall n \ge 1 $$
 which satisfies a set of ``good'' properties, stated in Proposition \ref{sec6}.2 (cf
\cite{Jo}, \cite{Jones}, \cite{IsOg2}, \cite{Co2}). Let's point out here on
another important fact that the Jucys--Murphy element $y_n(x)$ satisfies the reflection equation (\ref{reflH}). We also introduce a family of mutually commuting elements $\tau_n(x;;\vec{z}_(n)) \in BMW_n(C)$, the so-called {\it dressing
$JM$-operators} which are an analogue of the Sklyanin transfer matrices
\cite{Skl}, and the coefficients in the expansion of  $\tau_n(x;;\vec{z}_{(n)})$ over  the variable $x$  (for the homogeneous case $z_i =1, \forall$ ) are the Hamiltonians for the open Birman--Murakami--Wenzl chain models with nontrivial boundary conditions, see e.g. \cite{IsOg2}, and example at the end of Section \ref{sec6}.1. ~~\underline{Section \ref{sec6}.2} is devoted to construction of the Bethe subalgebras in the affine $BMW_n(C)$ algebras and a factorizibility property of the
corresponding $qKZ$ connections. We will show that the flat connections
${\sf A}_{i}'(z)$, see (\ref{Bax05}), are images under the map (\ref{Zam10K}) of certain
 elements ${\sf J}_i \in B_n(C)$ which under the special limit (\ref{Zam77}) one can deduce the $BMW$ analog (\ref{qkz03}) of the Cherednik's connections have been
introduced in \cite{Ch} for Hecke algebras.  As an application, in Section \ref{sec6} we construct a {\it baxterization}
of the type $A$ Dunkl--Cherednik elements $Y_i \in DAHA$, which have been
in-depth studied in \cite{Ch}.

\section{Affine braid groups of type
$A^{(1)},B^{(1)},C^{(1)},D^{(1)}$ and  Jucys--Murphy elements\label{ABG}}
\setcounter{equation}0

First consider affine braid group $B_n(C^{(1)})$
 with generators $\{ T_0, \dots , T_n \}$
subject to defining relations
 \be
 \lb{Affbg}
 T_i \, T_{i+1} \, T_i = T_{i+1} \, T_{i} \, T_{i+1} \; , \;\;\;\; i = 1, \dots, n-2 \; ,
 \ee
 \be
 \lb{Affbg2}
 \begin{array}{c}
 T_1 \, T_0 \, T_1 \, T_0 = T_0 \, T_1 \, T_0 \, T_1  \; , \;\;\;
 \\ [0.2cm]
 T_{n-1} \, T_n \, T_{n-1} \, T_n = T_n \, T_{n-1} \, T_n \, T_{n-1}  \; ,
 \end{array}
 \ee
 where $T_0,T_n$ --- two affine generators. Let $||m_{ij}||$ be symmetric matrix with
 integer coefficients $m_{ij} \geq 2$.
 The structure relations (\ref{Affbg}), (\ref{Affbg2}) of the group $B_n(C^{(1)})$ can be written as
 $\underbrace{T_i \, T_{j} \, T_i \cdots}_{m_{ij}}
  = \underbrace{T_{j} \, T_{i} \, T_{j} \cdots}_{m_{ji}}$
 and correspond to the Coxeter graph of the type $C^{(1)}$

 \unitlength=5mm
\begin{picture}(17,4.5)

\put(2,2){\circle{0.4}}
\put(1.5,2.5){$T_0$}
\put(4.5,2){\circle{0.4}}
\put(4.5,2.5){$T_1$}
\put(2.2,1.9){\line(1,0){2.1}}
\put(2.2,2.1){\line(1,0){2.1}}
\put(4.7,2){\line(1,0){2}}
\put(7,2){$. \; . \; . \; . \; . \; . \; .$}

\put(10.5,2){\circle{0.4}}
\put(10.5,2.5){$T_{n-2}$}
\put(10.7,2){\line(1,0){2}}
\put(13,2){\circle{0.4}}
\put(12.7,2.5){$T_{n-1}$}
\put(15.5,2){\circle{0.4}}
\put(15.5,2.5){$T_n$}
\put(13.2,1.9){\line(1,0){2.1}}
\put(13.2,2.1){\line(1,0){2.1}}

\end{picture}

\vspace{-1.5cm}

 \be
\lb{C1}
{}
\ee

\vspace{0.5cm}

\noindent
where the number of lines between nodes $i$ and $j$ is equal to $(m_{ij}-2)$. Note that for the
group $B_n(C^{(1)})$ defined by (\ref{Affbg}), (\ref{Affbg2}) we have two automorphisms 
$\rho_1$ and $\rho_2$:
 \be
\lb{autom}
\rho_1(T_i) = T_i^{-1} \; , \;\;\;\; \rho_2(T_i) = T_{n-i} \; .
\ee

The well known statement is:

{\bf Proposition \ref{ABG}.1.} {\it The affine braid group $B_n(C^{(1)})$ contains the
commutative subgroups which are generated by the following sets of elements
\be
\lb{jucys1}
\begin{array}{l}
\bullet~~J_i=
 \biggl( T_{i-1}^{-1} \cdots T_{1}^{-1}\biggr)~
 \biggl(T_0 \cdots T_n\biggr)~
 \biggl(T_{n-1} \cdots T_{i}\biggr)
  \; ,  \;\;\;\; i=1,\dots,n \; ,
 \\ [0.2cm]
 \bullet~~~~\overline{J}_i=
 \biggl(T_{i-1} \cdots T_{1}\biggr)~\biggl(T_0 \cdots T_{n}\biggr)~
 \biggl(T_{n-1}^{-1} \cdots T_{i}^{-1}\biggr)
\; , ~~i=1, \ldots, n \; ,
 \\ [0.2cm]
\bullet~~(Jucys-Murphy ~ elements)~~~
a_{i}:=   \biggl(T_{i-1} \cdots T_{1}\biggr)~T_0~\biggl(T_{1} \cdots T_{i-1}\biggr) \, ,
\;\;\;\; i=1,\dots,n \; , \\ [0.2cm]
\bullet~~(Jucys-Murphy ~ elements)~~~b_{i} := \biggl(T_{i} \cdots T_{n-1}\biggr)
 ~T_n~ \biggl(T_{n-1} \cdots T_{i}\biggr)
 ~~i=1, \ldots, n \; .
\end{array}
\ee}

{\bf Proof.} The proof of commutativity
of the elements $a_i$ is straightforward and follows from the fact that
$[a_i, \, T_j]=0$ for $i > j$. The commutativity of the elements $b_i$
follows from the commutativity of elements $a_i$ since we have $b_{n-i+1} = \rho_2(a_i)$,
where the automorphism $\rho_2$ is defined in (\ref{autom}).

 Now  we prove the commutativity of the elements $J_i$
 (it will be important for our consideration below). We introduce the element
  \be
 \lb{XX01}
 X = \prod\limits_{k=0}^{n}~T_k = T_0 \cdots T_{n} \; .
 \ee
 For this element we have the following identities
 \be
 \lb{XX02}
 \begin{array}{c}
 X \, T_i = T_{i+1} \, X \; , \;\;\; (i=1,\dots , n-2) \; , \\ [0.2cm]
 T_1 \cdot X^2 =  T_1 \cdot T_0 \, T_1 \, T_0 \, (T_2 \, T_1) \, (T_3 \, T_2)
 \cdots (T_{n-1} \, T_{n-2}) T_n \, T_{n-1} \, T_n =
 X^2 \,  T_{n-1} \; ,
 \end{array}
 \ee
 where in the proof of these identities we have used (\ref{Affbg}), (\ref{Affbg2}).
 With the help of the operator $X$ (\ref{XX01}) the element $\overline{J}_k$
 (\ref{jucys1}) can be written as
 $$
  \overline{J}_k =
  T_{k-1} \cdots T_1 \cdot X
 \cdot  T^{-1}_{n-1}  \cdots T_k^{-1} =
  T_{k-1} \cdots T_2 \cdot X
 \cdot  T_{n}^{-1} T^{-1}_{n-1}  \cdots T_k^{-1} \; .
 $$
Let $k > r$. Then by using (\ref{Affbg}), (\ref{Affbg2}) and (\ref{XX02})
 we have
 $$
  \begin{array}{c}
 \overline{J}_k \;  \overline{J}_r =  (T_{k-1} \cdots T_1 \cdot X
 \cdot  T^{-1}_{n-1}  \cdots T_k^{-1}) \cdot
 (T_{r-1} \cdots T_1 \cdot X
 \cdot  T^{-1}_{n-1}  \cdots T_r^{-1}) = \\ [0.2cm]
  =  (T_{k-1} \cdots T_1) \cdot X
 \cdot  (T_{r-1} \cdots T_1) \cdot (T^{-1}_{n-1}  \cdots T_k^{-1}) \cdot
  X \cdot  (T^{-1}_{n-1}  \cdots T_r^{-1}) = \\ [0.2cm]
  =  (T_{k-1} \cdots T_1) \cdot  (T_{r} \cdots T_2) \cdot X
  \cdot  X \cdot
 (T^{-1}_{n-2}  \cdots T_{k-1}^{-1}) \cdot (T^{-1}_{n-1}  \cdots T_r^{-1}) =
 \\ [0.2cm]
  = (T_{r-1} \cdots T_1) \cdot (T_{k-1} \cdots T_1) \cdot X^2
   \cdot (T^{-1}_{n-1}  \cdots T_r^{-1}) \cdot
 (T^{-1}_{n-1}  \cdots T_{k}^{-1}) = \\ [0.2cm]
  = (T_{r-1} \cdots T_1) \cdot (T_{k-1} \cdots T_2) \cdot X^2
  \cdot  T_{n-1}  \cdot (T^{-1}_{n-1}  \cdots T_r^{-1}) \cdot
 (T^{-1}_{n-1}  \cdots T_{k}^{-1}) = \\ [0.2cm]
  = (T_{r-1} \cdots T_1)  \cdot X \cdot (T_{k-2} \cdots T_1)
  \cdot (T^{-1}_{n-1}  \cdots T_{r+1}^{-1})  \cdot X
  \cdot (T^{-1}_{n-1}  \cdots T_{k}^{-1}) =
  \\ [0.2cm]
  = (T_{r-1} \cdots T_1)  \cdot X
  \cdot (T^{-1}_{n-1}  \cdots T_{r}^{-1}) \cdot (T_{k-1} \cdots T_1)
    \cdot X \cdot (T^{-1}_{n-1}  \cdots T_{k}^{-1})
     = \overline{J}_r \;  \overline{J}_k  \; ,
  \end{array}
 $$
 where to obtain the last line we use the identity $(k > r)$
 $$
  \begin{array}{c}
 (T_{k-2} \cdots T_1)
  \cdot (T^{-1}_{n-1}  \cdots T_{r+1}^{-1}) =
  (T_{k-2} \cdots T_r) \cdot  (T_{r-1} \cdots T_1)  \cdot
  (T^{-1}_{n-1}  \cdots T_{k}^{-1}) (T^{-1}_{k-1}  \cdots T_{r+1}^{-1})
   = \\ [0.2cm]
  = (T^{-1}_{n-1}  \cdots T_{k}^{-1}) (T_{k-2} \cdots T_r)
  \cdot (T^{-1}_{k-1}  \cdots T_{r+1}^{-1}) (T_{r-1} \cdots T_1) = \\ [0.2cm]
  = (T^{-1}_{n-1}  \cdots T_{k}^{-1} T_{k-1}^{-1})
  (T_{k-1} T_{k-2} \cdots T_r)
  \cdot (T^{-1}_{k-1}  \cdots T_{r+1}^{-1}) (T_{r-1} \cdots T_1) = \\ [0.2cm]
   = (T^{-1}_{n-1}  \cdots  T_{k-1}^{-1})
   (T^{-1}_{k-2}  \cdots T_{r}^{-1})  \cdot
  (T_{k-1}  \cdots T_r)
  (T_{r-1} \cdots T_1) =
  (T^{-1}_{n-1}  \cdots T_{r}^{-1}) \cdot (T_{k-1} \cdots T_1) \; .
  \end{array}
 $$

  The commutativity of the elements $J_i$ follows from
  the commutativity of the elements $\overline{J}_i$ since we have
  $\rho_1(\rho_2(\overline{J}_{n-i+1})) = J_i^{-1}$, where automorphisms
  $\rho_1$ and $\rho_2$ are defined in (\ref{autom}). \hfill \qed

\vspace{0.2cm}

The quotient of the group $B_n(C^{(1)})$ by
additional relations $T_i^2 =1$ $(\forall i)$ is called Coxeter group of
the type $C^{(1)}$. This group is denoted as $W_n(C^{(1)})$.
At the end of this Section we present the explicit realization of $W_n(C^{(1)})$
 which we use below.
 Introduce the set of spectral parameters $(z_1 , \dots , z_n)$, $z_i \in \mathbb{C}$.
 Now we define a representation ${\sf s}: \; T_i \to s_i$ of $B_{n}$:
  \be
 \lb{Affbg1}
 \begin{array}{c}
 s_i \; : \; (z_1, \dots, z_i, z_{i+1} , \dots , z_n) \;\; \to \;\;
 (z_1, \dots, z_{i+1}, z_{i} , \dots , z_n) \;\;\;\;\; (i=1,\dots,n-1) \; ,   \\ [0.2cm]
 s_0  \; : \; (z_1,z_2, \dots , z_n) \;\; \to \;\;
  (\sigma(z_1),z_2, \dots, z_n)  \; , \\ [0.2cm]
 s_n  \; : \; (z_1, \dots , z_{n-1} , z_n) \;\; \to \;\;
  (z_1, \dots , z_{n-1} , \bar{\sigma}(z_n) )  \; ,
 \end{array}
 \ee
 where $\sigma$, $\bar{\sigma}$ are two involutive mappings $\mathbb{C} \to \mathbb{C}$
 such that
 $(\sigma)^2=1$, $(\bar{\sigma})^2 =1$. We specify these involutions in next
 Sections. From (\ref{Affbg1}) one can check
 that operators $s_0,s_i,s_n$ satisfy (\ref{Affbg}), (\ref{Affbg2}) and
 moreover we have $s_0^2 = s_n^2 =s_i^2=1$.
 Thus, equations (\ref{Affbg1}) give the representation of the Coxeter group $W_n(C^{(1)})$. For special choices of
 $\sigma$ and $\bar{\sigma}$,
 namely $\sigma(z) = 1-z$ and $\bar{\sigma}(z) =-z$, the
 representation (\ref{Affbg1}) have been used in \cite{Ch},\cite{Stok}.

   \vspace{0.2cm}

  \noindent
  {\bf Remark 1.} Denote by $B_n(C)$ the subgroup of the affine
  braid group $B_n(C^{(1)})$
  generated by elements $T_i$ $(i=0,\dots, n-1)$ with defining relations given in (\ref{Affbg})
  and in first line of (\ref{Affbg2}). The group $B_n(C)$ is associated to
  the Coxeter graph of $C$-type

   \unitlength=5mm
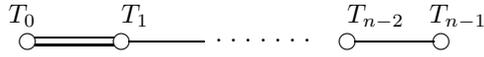
\begin{picture}(17,4.5)

\put(2,2){\circle{0.4}}
\put(1.5,2.5){$T_0$}
\put(4.5,2){\circle{0.4}}
\put(4.5,2.5){$T_1$}
\put(2.2,1.9){\line(1,0){2.1}}
\put(2.2,2.1){\line(1,0){2.1}}
\put(4.7,2){\line(1,0){2}}
\put(7,2){$. \; . \; . \; . \; . \; . \; .$}

\put(10.5,2){\circle{0.4}}
\put(10.5,2.5){$T_{n-2}$}
\put(10.7,2){\line(1,0){2.1}}
\put(13,2){\circle{0.4}}
\put(12.7,2.5){$T_{n-1}$}

\end{picture}

 \noindent
  Consider the homomorphism (projection)
  $\rho$: $B_n(C^{(1)}) \to B_n(C)$
  such that $\rho(T_i) = T_i$ $(i=0,\dots, n-1)$ and $\rho(T_n) =1$.
 It is clear that under this projection
  we have $a_i = \rho(\overline{J}_i)$
  and it means that the commutativity of $a_i$
   follows from the commutativity of $\overline{J}_i$.
   The elements $a_i$ given in (\ref{jucys1}) generate the commutative set
  in the subgroup $B_n(C) \subset B_n(C^{(1)})$.

     \vspace{0.2cm}

  \noindent
  {\bf Remark 2.} Denote by $B_n(A^{(1)})$ the affine braid group
  which corresponds to the affine $A$-type Coxeter graph

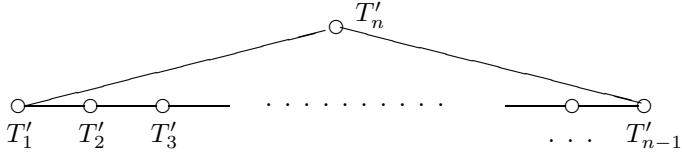
\unitlength=8mm
\begin{picture}(17,3)
\put(1.9,1.1){\circle{0.2}}
\put(1.7,0.5){$T'_1$}
\put(2,1.1){\line(1,0){1}}
\put(3.1,1.1){\circle{0.2}}
\put(2.9,0.5){$T'_2$}
\put(3.2,1.1){\line(1,0){1}}
\put(4.3,1.1){\circle{0.2}}
\put(4.1,0.5){$T'_3$}
\put(4.4,1.1){\line(1,0){1}}
\put(6,1.1){. . . . . . . . . .}
\put(10,1.1){\line(1,0){1}}
\put(11.1,1.1){\circle{0.2}}
\put(10.7,0.5){. . .}
\put(11.2,1.1){\line(1,0){1}}
\put(12.3,1.1){\circle{0.2}}
\put(12,0.5){$T'_{n-1}$}

\put(2,1.1){\line(4,1){5}}
\put(7.25,2.4){\line(4,-1){5}}
\put(7.18,2.4){\circle{0.2}}
\put(7.5,2.5){$T'_n$}
\end{picture}

\noindent
We call group $B_n(A^{(1)})$ $(n>2)$ a periodic $A$-type braid group.
This group is generated by invertible elements $T'_i$ $(i=1,\dots,n)$
and according to its Coxeter graph we have the defining relations
 \be
 \lb{perBn}
 T'_i \, T'_{i+1} \, T'_i = T'_{i+1} \, T'_{i} \, T'_{i+1} \; , \;\;\;\; i = 1, \dots, n \; , \\ [0.2cm]
 \ee
 where we impose the periodic conditions $T'_{i +n}=T'_i$.

 Note that the group $B_n(A^{(1)})$ possesses automorphisms
 \be
 \lb{perBn1}
 \rho_3( T'_i) = T'_{i+1} \; , \;\;\; \rho_4( T'_i) = T'_{n-i+1} \; , \;\;\; \rho_5( T'_i) = T_{i}^{\prime -1} \; .
 \ee
 Define the extension $\bar{B}_n(A^{(1)})$ of the group $B_n(A^{(1)})$
 by adding an additional generator $\bar{X}$  with defining relations (cf. (\ref{XX02}))
 \be
 \lb{perBn2}
 \bar{X} \, T'_i = T'_{i+1} \,\bar{X} \;\;\;\;\;
 (i=1,\dots, n) \;\;\;\; \Rightarrow \;\;\;\;
 T'_1 \cdot \bar{X}^2 = \bar{X}^2 \cdot T'_{n-1} \; .
 \ee
 Namely, we add operator $\bar{X}$ which serves the
 automorphism $\rho_3$:
 $\rho_3(T'_i) = \bar{X} \, T'_i \bar{X}^{-1}$ in (\ref{perBn1}).
 Then for the group $\bar{B}_n(A^{(1)})$ one can construct the following
 commuting sets of elements
 \be
 \lb{perBn3}
 \begin{array}{c}
 J'_k = T^{\prime -1}_{k-1} \cdots T^{\prime -1}_1
 \cdot \bar{X}
 \cdot  T'_{n-1}  \cdots T'_k \;\;\;\;\; (k=1,\dots,n) \; , \\ [0.2cm]
  \bar{J}'_k = \rho_5(J'_k) = T'_{k-1} \cdots T'_1 \cdot \bar{X}
 \cdot  T^{\prime -1}_{n-1}  \cdots T_k^{\prime -1} \;\;\;\;\; (k=1,\dots,n) \; ,
 \end{array}
 \ee
 where we have defined $\rho_5(\bar{X}) = \bar{X}$ (this is compatible with (\ref{perBn2})).

 Now we introduce the element $\bar{T}_n$ in $B_n(C^{(1)})$ as following
 \be
 \lb{XX03}
 \bar{T}_n :=
 X^{-1} \, T_1 \cdot X = X \,  T_{n-1} \, X^{-1} \;\; \in \;\; B_n(C^{(1)}) \; ,
 \ee
 where $X$ is given in (\ref{XX01}).
 The element (\ref{XX03}) satisfies
 periodic braid relations
 $$
 \bar{T}_n \,  T_{n-1} \,  \bar{T}_n = T_{n-1} \, \bar{T}_n \, T_{n-1} \; , \;\;\;
 \bar{T}_n \,  T_{1} \,  \bar{T}_n = T_{1} \, \bar{T}_n \, T_{1} \; ,
 $$
 where we have used (\ref{XX02}). Thus, we have the homomorphic maps
  (embeddings) $\rho'$:
 $B_n(A^{(1)}) \to B_n(C^{(1)})$ and
 $\rho^{\prime \prime}$: $\bar{B}_n(A^{(1)}) \to B_n(C^{(1)})$ such that
 $$
  \begin{array}{c}
 \rho'(T'_i) = T_i \;\;\;\; (i=1,\dots,n-1) \; , \;\;\;\;\;
 \rho'(T'_n) = \bar{T}_n \; , \\ [0.2cm]
  \rho^{\prime \prime}(T'_i) = T_i \;\;\;\;
  (i=1,\dots,n-1) \; , \;\;\;\;\;
 \rho^{\prime \prime}(T'_n) = \bar{T}_n
 \; , \;\;\; \rho^{\prime \prime}(\bar{X}) = X \; .
  \end{array}
 $$
 It means that
 $B_n(A^{(1)})$ and
$\bar{B}_n(A^{(1)})$ are subgroups in $B_n(C^{(1)})$
with generators $(T_1,\dots,T_{n-1},\bar{T}_n)$
and $(T_1,\dots,T_{n-1},\bar{T}_n, X)$, respectively.

 \vspace{0.2cm}

  \noindent
  {\bf Remark 3.} Consider the braid group $B_{n+1}(B^{(1)})$ which is associated to the graph

     \unitlength=5mm
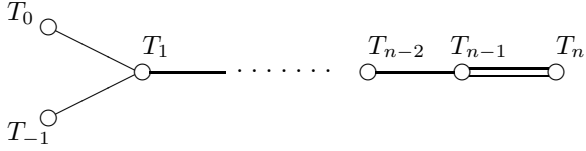
\begin{picture}(17,4.5)

\put(2,3.2){\circle{0.4}}
\put(2,0.8){\circle{0.4}}
\put(0.9,3.4){$T_0$}
\put(0.9,0.2){$T_{-1}$}
\put(4.5,2){\circle{0.4}}
\put(4.5,2.5){$T_1$}
\put(2.2,3.1){\line(2,-1){2.1}}
\put(2.2,0.9){\line(2,1){2.1}}
\put(4.7,2){\line(1,0){2}}
\put(7,2){$. \; . \; . \; . \; . \; . \; .$}

\put(10.5,2){\circle{0.4}}
\put(10.5,2.5){$T_{n-2}$}
\put(10.7,2){\line(1,0){2.1}}
\put(13,2){\circle{0.4}}
\put(12.7,2.5){$T_{n-1}$}
\put(15.5,2){\circle{0.4}}
\put(15.5,2.5){$T_n$}
\put(13.2,1.9){\line(1,0){2.1}}
\put(13.2,2.1){\line(1,0){2.1}}

\end{picture}

\vspace{0.2cm}

 \noindent
The defining relations for this group are
 \be
 \lb{Affbd}
 \begin{array}{c}
  T_i \, T_{i+1} \, T_i = T_{i+1} \, T_{i} \, T_{i+1} \; , \;\;\;\; i = 0,1, \dots, n-1 \; ,
  \\ [0.2cm]
 T_{-1} \, T_1 \, T_{-1} =  T_1 \, T_{-1} \, T_1  \;\; , \;\;\;\;\;\;
 T_{-1} \, T_0 =  T_0 \, T_{-1}   \; ,
 \\ [0.2cm]
 T_{n-1} \, T_n \, T_{n-1} \, T_n = T_n \, T_{n-1} \, T_n \, T_{n-1}  \; .
 \end{array}
 \ee
 Introduce the element
 \be
 \lb{tilT0}
 \tilde{T}_0 = T_{-1} \, T_0 \; ,
 \ee
 which
 in view of (\ref{Affbd}) satisfies relation
 \be
 \lb{Affbf}
 \tilde{T}_0 \, T_1 \, \tilde{T}_0 \, T_1 = T_1 \, \tilde{T}_0 \, T_1 \, \tilde{T}_0
 \ee
 So, $B_{n}(C^{(1)})$ is a subgroup in
 $B_{n+1}(B^{(1)})$ and we have the homomorphism
 (embedding)  $\tilde{\rho}$: $B_{n}(C^{(1)})\to B_{n+1}(B^{(1)})$ which is defined by
 the map
 \be
 \lb{mapT}
 \tilde{\rho} \; : \;\;
 T_0 \; \to \; \tilde{T}_0 \; , \;\;\;\; T_i \; \to \; T_i
 \;\;\;\;\;\; (i=1,\dots,n) \; .
 \ee
 Thus, according to the Proposition \ref{ABG}.1 we have the following commuting sets
 for the group $B_{n+1}(B^{(1)})$
 \be
\lb{jucys11}
\begin{array}{l}
 \tilde{J}_i= \biggl(\prod\limits_{k=i-1}^1 \; T_{k}^{-1}  \biggr)~
 \tilde{X}~
\biggl( \prod\limits_{k=n-1}^i \; T_{k}  \biggr)
  \;   \;\;\;\; (i=1,\dots,n) \; ,
 \\ [0.2cm]
 \overline{\tilde{J}}_i= \biggl( \prod\limits_{k=i-1}^1 \; T_{k} \biggr)~\tilde{X}~
\biggl( \prod\limits_{k=n-1}^i \; T_{k}^{-1}  \biggr)
\;~~(i=1, \ldots, n) \; ,
\end{array}
\ee
where $\tilde{X}=\tilde{T}_0 \, T_1 \cdots T_n$
is the image of the element $X \in B_{n+1}(C^{(1)})$
presented in (\ref{XX01}).

 \vspace{0.2cm}

  \noindent
  {\bf Remark 4.} The braid group $B_{n+2}(D^{(1)})$ which is associated with the graph

     \unitlength=5mm
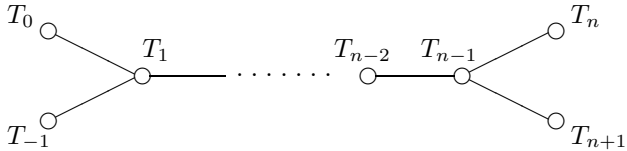
\begin{picture}(17,4.5)

\put(2,3.2){\circle{0.4}}
\put(2,0.8){\circle{0.4}}
\put(0.9,3.4){$T_0$}
\put(0.9,0.2){$T_{-1}$}
\put(4.5,2){\circle{0.4}}
\put(4.5,2.5){$T_1$}
\put(2.2,3.1){\line(2,-1){2.1}}
\put(2.2,0.9){\line(2,1){2.1}}
\put(4.7,2){\line(1,0){2}}
\put(7,2){$. \; . \; . \; . \; . \; . \; .$}

\put(10.5,2){\circle{0.4}}
\put(9.6,2.5){$T_{n-2}$}
\put(10.7,2){\line(1,0){2.1}}
\put(13,2){\circle{0.4}}
\put(11.9,2.5){$T_{n-1}$}
\put(13.2,2.1){\line(2,1){2.1}}
\put(13.2,1.9){\line(2,-1){2.1}}
\put(15.9,3.4){$T_n$}
\put(15.9,0.2){$T_{n+1}$}
\put(15.5,3.2){\circle{0.4}}
\put(15.5,0.8){\circle{0.4}}

\end{picture}

\vspace{0.2cm}

 \noindent
has defining relations
 \be
 \lb{Affdd}
 \begin{array}{c}
  T_i \, T_{i+1} \, T_i = T_{i+1} \, T_{i} \, T_{i+1} \; , \;\;\;\; i = 0,1, \dots, n \; ,
  \\ [0.2cm]
 T_{-1} \, T_1 \, T_{-1} =  T_1 \, T_{-1} \, T_1  \;\; , \;\;\;\;\;\;
 T_{-1} \, T_0 =  T_0 \, T_{-1}   \; ,
 \\ [0.2cm]
  T_{n-1} \, T_{n+1} \, T_{n-1} =
  T_{n+1} \, T_{n-1} \, T_{n+1}  \;\; , \;\;\;\;\;\;
 T_{n} \, T_{n+1} =  T_{n+1} \, T_{n}   \; .
 \end{array}
 \ee
 Note that the element $\tilde{T}_n = T_n \, T_{n+1}$
 obeys relations
 $$
 \tilde{T}_n \, T_{n-1} \, \tilde{T}_n \, T_{n-1} =
 T_{n-1} \, \tilde{T}_n \, T_{n-1} \, \tilde{T}_n \; .
 $$
 Thus the elements
 $(T_{-1},T_0,T_1,\dots,T_{n-1},\tilde{T}_n)$ generate
 the subgroup $B_{n+1}(B^{(1)})$ in $B_{n+2}(D^{(1)})$
 and we have the homomorphism
 (embedding) $\rho_0:$ $B_{n+1}(B^{(1)}) \to B_{n+2}(D^{(1)})$
 such that
 \be
 \lb{mapT2}
 \rho_0: T_i \;\; \to \;\; T_i \;\;\;\;\;
 (i=-1,0,1,\dots,n-1) \; , \;\;\;
 \rho_0: T_n \;\; \to \;\; \tilde{T}_n \; .
 \ee
Define the element (cf. (\ref{XX01}))
$$
X^{\prime \prime} = \tilde{T}_0 \, T_1 \cdots T_{n-1} \,
\tilde{T}_n \; ,
$$
where $\tilde{T}_0$ is defined as in (\ref{tilT0}).
Then we again have two sets of commuting elements
(cf. (\ref{jucys1}), (\ref{jucys11}))
  \be
\lb{jucys5}
\begin{array}{l}
 J_i^{\prime \prime}
 = \biggl(\prod\limits_{k=i-1}^1 \; T_{k}^{-1}  \biggr)~
 X^{\prime \prime}~
\biggl( \prod\limits_{k=n-1}^i \; T_{k}  \biggr)
  \;   \;\;\;\; (i=1,\dots,n) \; ,
 \\ [0.2cm]
 \overline{J}_i^{\prime \prime}=
 \biggl( \prod\limits_{k=i-1}^1 \; T_{k} \biggr)~
 X^{\prime \prime}~
\biggl( \prod\limits_{k=n-1}^i \; T_{k}^{-1}  \biggr)
\;~~(i=1, \ldots, n) \; .
\end{array}
\ee

 Finally we stress that the quotient of the group $B_{n+2}(D^{(1)})$ with respect to
 the relations $T_0 = T_{-1}$ (or $T_n = T_{n+1}$) is
 isomorphic to the braid group $B_{n+2}(D)$
 associated to the Coxeter graph of classical $D$-type.
 The commutative elements in this case are given by the same formulas as in (\ref{jucys11}), where instead of $\tilde{X}$
 we have to substitute
 element $X(D) = T_0^2 T_1 \cdots T_{n-1} \tilde{T}_n$
 (or $X(D) = \tilde{T}_0 \, T_1 \cdots T_{n-1} T_n^2$).

 \vspace{1cm}

\section{General picture}
\setcounter{equation}0

 \paragraph{1. Affine root systems and affine
  Weyl groups (see \cite[Section~1]{C1}).}


 ${}$

Let $R_n$ be a root system of type $A_n,B_n,\allowbreak\ldots,\allowbreak F_n,G_n$. We will write $R$ also
for the type of the root system. Let $\alpha_1,\allowbreak\ldots,\allowbreak\alpha_n\in R_n$ be simple roots,
$\omega^\vee_1,\allowbreak\ldots,\allowbreak \omega^\vee_n$ \!--- fundamental coweights, $(\omega^\vee_i,\alpha_i)=\delta_i^j$,
\,$\theta$ --- the maximal root. The Dynkin diagram of the affine root system
$R_n^{(1)}\!$ is obtained by adding the root $-\,\theta$ to the simple roots
$\alpha_1,\allowbreak\ldots,\allowbreak\alpha_n$. The affine simple root is $\alpha_0=[-\,\theta,1]$ in the notation
of \cite{C1}.

\smallskip
For $\alpha\in R_n$, denote $\alpha^\vee=2\alpha/(\alpha,\alpha)$. Let
$Q^\vee\!=\bigoplus_{i=1}^n\mathbb{Z}\alpha^\vee_i$ be the coroot lattice,
$P^\vee\!=\bigoplus_{i=1}^n \mathbb{Z}\,\omega^\vee_i$ be the coweight lattice, and
$P^\vee_+\!=\bigoplus_{i=1}^n \mathbb{Z}_{\ge0}\omega^\vee_i$

\smallskip
Let $s_\alpha$ be the reflection corresponding to a root $\alpha\in R_n^{(1)}\!$,
and $s_i=s_{\alpha_i}$.
The Weyl group $W$ of type $R_n$ is generated by the reflections $s_1,\allowbreak\ldots,\allowbreak s_n$.

\smallskip
The affine Weyl group $W^{(a)}$ of type $R_n^{(1)}\!$ is generated by
the reflections $s_0,s_1,\allowbreak\ldots,\allowbreak s_n$ and is isomorphic to the semidirect product
$W\ltimes Q^\vee$, with $s_0=\theta^\vee s_\theta$. Here we identify $W$ and
$Q^\vee\!$ with the respective subgroups of $W^a\!$.

\smallskip
The extended affine Weyl group $W^{(b)}$ of type $R_n^{(1)}\!$ is the semidirect
product $\,\widetilde{W}=W\ltimes P^\vee$. It is also isomorphic to the semidirect
product $\Pi\ltimes W^a$, where $\Pi=P^\vee/Q^\vee$. The elements of the
subgroup $\Pi\subset\widetilde{W}$ ``permute'' the reflections $s_0,\allowbreak\ldots,\allowbreak s_n$ \!--- for any
$i$ and $\pi\in\Pi$, $\,\pi s_i\pi^{-1}=s_j$ \,for some $j=\pi[i]$.

\smallskip
Define the length on $\widetilde{W}$ by $\ell(s_i)=1$ and $\ell(\pi)=0$ for
$\pi\in\Pi$. Then for $b,b'\!\in P^\vee_+\!\subset\widetilde{W}$,
\be
\label{ellb}
\ell(b+b')\,=\,\ell(b)+\ell(b')\,,
\ee
see \cite[Proposition~1.4]{C1}.

\smallskip
The affine braid group $B(R_n^{(1)})$ is generated by the elements
$S_0,\allowbreak\ldots,\allowbreak S_n$ subject to the same braid relations as $s_0,\allowbreak\ldots,\allowbreak s_n$
(we use $S_i$ to keep distinction from the generators $T_i$ in Section \ref{ABG}.)
The extended affine braid group $\widetilde{B}(R_n^{(1)})$ is the semidirect product
$\Pi\ltimes B(R_n^{(1)})$ --- for any $i$ and $\pi\in\Pi$,
$\,\pi S_i\pi^{-1}=S_{\pi[i]}$, (cf.~with relations \,(i), \!(ii)
in \cite[Definition~3.1]{C1}).

\smallskip
For $\widetilde{w}\in\widetilde{W}$ with a reduced decomposition
$\widetilde{w}=\pi s_{i_1}\dots s_{i_k}$, $\pi\in\Pi$, $k=\ell(\widetilde{w})$,
the element $S_{\widetilde{w}}=\pi S_{i_1}\dots S_{i_k}\in\widetilde{W}$ does not depend
on the reduced decomposition, and $S_{\widetilde{w}\widetilde{w}'}=S_{\widetilde{w}}S_{\widetilde{w}'}$ provided
$\ell(\widetilde{w}\widetilde{w}')=\ell(\widetilde{w})+\ell(\widetilde{w}')$, $\,\widetilde{w},\widetilde{w}'\in\widetilde{W}$. Hence,
the elements $S_b$, $b\in P^\vee_+\!\subset\widetilde{W}$ generate a commutative subgroup
of $\widetilde{W}$ because $S_bS_{b'}=S_{b+b'}=S_{b'}S_b$ for any
$b,b'\!\in P^\vee_+\!\subset\widetilde{W}$, see \eqref{ellb}.
(Cf.~with \cite[formula (3.8)]{C1}.)

\smallskip
For fundamental coweights $\omega^\vee_1,\allowbreak\ldots,\allowbreak \omega^\vee_n$, set
\be
\label{Yi}
Y_i\,=\,S_{\omega^\vee_i},\qquad i=1,\allowbreak\ldots,\allowbreak n\,.
\ee
The elements \,$Y_1,\allowbreak\ldots,\allowbreak Y_n\in\widetilde{W}$ pairwise commute.

\medskip
\paragraph{2. Groups $\widehat{B}(C_n^{(1)})$ and $\widetilde{B}(C_n^{(1)})$.}

${}$

The group $\widehat{B}(C_n^{(1)})$ is generated by the elements
$T_0,\allowbreak\ldots,\allowbreak T_n\in B(C_n^{(1)})$, see~(\ref{Affbg}),
 (\ref{Affbg2}) and by the element $U$ with relations
\be
\label{UT}
UT_iU^{-1}=\,T_{n-i}\,,\quad i=0,\allowbreak\ldots,\allowbreak n\,.
\ee
In other words, $UGU^{-1}=\rho_2(G)$ for any $G\in B(C_n^{(1)})$, where
$\rho_2$ is given by formula~(\ref{autom}). The element $U^2$
is central.

\smallskip
Set $I_i\,=\,J_1\dots J_i$, \,$i=1,\allowbreak\ldots,\allowbreak n$, where $J_1,\allowbreak\ldots,\allowbreak J_n$ are given
by~(\ref{jucys1}). Also,
\be
\label{Ii}
I_i\,=\,(X T_{n-1}\dots T_i)^{\,i},\qquad i=1,\allowbreak\ldots,\allowbreak n\,,
\ee
where $X=T_0\dots T_n$, \,see~(\ref{XX01}). Let
\be
\label{Zn}
Z\,=\,T_0\dots T_{n-1}\;T_0\dots T_{n-2}\,\dots\,T_0\,T_1\;T_0\,U\,.
\ee
The element $Z$ commutes with $T_1,\allowbreak\ldots,\allowbreak T_{n-1}$ and $X$, and hence
by~\eqref{Ii}, commutes with $I_1,\allowbreak\ldots,\allowbreak I_n$. Moreover, $Z^2=I_n U^2$.
One more nice formula
\be
\label{Io}
I_i\,=\,
X^{i}\,T_{n-i}\dots T_1\;T_{n-i+1}\dots T_2\,\dots\,T_{n-1}\dots T_i\,.
\ee

\smallskip
The group $\widetilde{B}(C_n^{(1)})$ is the quotient of $\widehat{B}(C_n^{(1)})$ by relation
$U^2=1$. The identification is $S_i=T_i$, \,$i=0,\allowbreak\ldots,\allowbreak n$, and $\Pi=\{1,U\}$.
Also, $Y_i=I_i$, \,$i=1,\allowbreak\ldots,\allowbreak n-1$, and $Y_n=Z$.

\medskip
\paragraph{3. Groups $B(B_n^{(1)})$ and $\widetilde{B}(B_n^{(1)})$.}

${}$

The group $\widetilde{B}(B_n^{(1)})$ is the quotient of $B(C_n^{(1)})$ by relation
$T_0^2=1$. The identification is $S_i=T_i$, \,$i=1,\allowbreak\ldots,\allowbreak n$,
\,$S_0=T_0T_1T_0$, and $\Pi=\{1,T_0\}$. Thus $S_0S_1=S_1S_0$
and $S_0S_2S_0=S_2S_0S_2$. Also $Y_i=I_i$, \,$i=1,\allowbreak\ldots,\allowbreak n$.
The commutative subgroup in $B(B_n^{(1)})$ is generated by the products
$J_1J_i=I_1\,I_i\,I_{i-1}^{-1}$, \,$i=1,\allowbreak\ldots,\allowbreak n$. Here $I_0=1$.

\smallskip
The relation with elements~(\ref{jucys11}), (\ref{jucys5}) is explained farther.

\medskip
 \paragraph{4. Groups $B(D_n^{(1)})$ and $\widetilde{B}(D_n^{(1)})$.}

 ${}$

The groups $B(D_n^{(1)})$ and $\widetilde{B}(D_n^{(1)})$ are subquotients of
$\widetilde{B}(C_n^{(1)})$. Let $\widetilde{B}'(C_n^{(1)})$ be the quotient of $\widetilde{B}(C_n^{(1)})$ by
relations $T_0^2=1$, \,$T_n^2=1$. (Recall that $U^2=1$ in $\widetilde{B}(C_n^{(1)})$.)
The subgroup $B(D_n^{(1)})\subset\widetilde{B}'(C_n^{(1)})$ is generated by
$S_0=T_0T_1T_0$, \,$S_n=T_nT_{n-1}T_n$,
and $S_i=T_i$, \,$i=1,\allowbreak\ldots,\allowbreak n-1$.

\smallskip
Let $\Pi_n=\{1,T_0U,(T_0U)^2,(T_0U)^3\}=
\{1,T_0U,T_0T_n,T_nU\}$ if $n$ is odd, and
$\Pi_n=\{1,T_0T_n,U,T_0T_nU\}$ if $n$ is even.
The subgroup $\widetilde{B}(D_n^{(1)})\subset\widetilde{B}'(C_n^{(1)})$ is generated by
$B(D_n^{(1)})$ and $\Pi=\Pi_n$.

\smallskip
Also \,$Y_i=I_i$, \,$i=1,\allowbreak\ldots,\allowbreak n-2$, $\,Y_{n-1}=I_{n-1}Z^{-1}$ and
\,$Y_n=Z$. The subgroup $\Pi_n$ can be recovered from the requirement
that $I_1$ and $Z$ belong to the subgroup generated by $B(D_n^{(1)})$ and
$\Pi=\Pi_n$.

\smallskip
The commutative subgroup in $B(D_n^{(1)})$ is generated by the products
$J_1J_i=I_1\,I_i\,I_{i-1}^{-1}$, \,$i=1,\allowbreak\ldots,\allowbreak n$. Here $I_0=1$.

\medskip
 \paragraph{5. Recursive definition of $I_1,\allowbreak\ldots,\allowbreak I_n$.}

 ${}$

Let $T'_0,\allowbreak\ldots,\allowbreak T'_{n-1}$ denote the generators of $B(C_{n-1}^{(1)})$,
and similarly for $I'_1,\allowbreak\ldots,\allowbreak I'_{n-1}$. There is an embedding
\be
\label{mu}
\mu:B(C_{n-1}^{(1)})\,\to\,B(C_n^{(1)})\,,
\ee
$$
\mu(T'_i)\,=\,T_i\,,\quad i=0,\allowbreak\ldots,\allowbreak n-2\,,\qquad
\mu(T'_{n-1})\,=\,T_{n-1}T_nT_{n-1}\,.
$$
Then
\be
\label{muI}
\mu(I'_i)\,=\,I_i\,,\quad i=1,\allowbreak\ldots,\allowbreak n-1\,.
\ee
This suggests a proof that the elements $I_1,\allowbreak\ldots,\allowbreak I_{n-1}, Z$ pairwise commute.
Since $Z^2=I_n U^2$, by induction it suffices to prove only that $Z$
commutes with $I_1,\allowbreak\ldots,\allowbreak I_{n-1}$. This follows from the fact that $Z$ commutes
with $T_1,\allowbreak\ldots,\allowbreak T_{n-1}$ and $X$, and formula \eqref{Ii}.



\smallskip
The recursive definition of $I_1,\allowbreak\ldots,\allowbreak I_n$ is reminiscent of the construction of Gelfand-Zetlin subalgebras, though $I_n$ is not central.


\medskip
 \paragraph{6. Relation with elements~(\ref{jucys11}), (\ref{jucys5}).}

 ${}$

To get elements~(\ref{jucys11}) compose the embedding $\mu$ with
automorphisms~(\ref{autom}) for $B(C_{n-1}^{(1)})$ and $B(C_n^{(1)})$:
$\lambda=\rho_2\circ\mu\circ\rho'_2$,
\be
\label{la}
\lambda(T'_i)\,=\,T_{i+1}\,,\quad i=1,\allowbreak\ldots,\allowbreak n-1\,,\qquad
\lambda(T'_0)\,=\,T_1T_0T_1\,.
\ee
The elements of the image $\lambda\bigl(B(C_{n-1}^{(1)})\bigr)$ commute with
$T_0$. The next formulae define one more embedding
$\tilde\lambda:B(C_{n-1}^{(1)})\to B(C_n^{(1)})$:
\be
\label{lat}
\tilde\lambda(T'_i)\,=\,T_{i+1}\,,\quad i=1,\allowbreak\ldots,\allowbreak n-1\,,\qquad
\tilde\lambda(T'_0)\,=\,T_0T_1T_0T_1\,.
\ee
Taking the quotient by the relation $T_0^2=1$ projects $B(C_n^{(1)})$ into
$\widetilde{B}(B_n^{(1)})$ and formulae \eqref{lat} into
\be
\label{latb}
\tilde\lambda(T'_i)\,=\,S_{i+1}\,,\quad i=1,\allowbreak\ldots,\allowbreak n-1\,,\qquad
\tilde\lambda(T'_0)\,=\,S_0S_1\,.
\ee
(Recall that \,$S_0=T_0T_1T_0$ and \,$S_i=T_i$, \,$i=1,\allowbreak\ldots,\allowbreak n$.)
Formulae \eqref{latb} coincide with the embedding $\tilde\rho$ in (\ref{mapT})
up to relabeling of generators.

\smallskip
To get elements~(\ref{jucys5}), the game is similar. First take an embedding
$B(C_{n-2}^{(1)})\to B(C_n^{(1)})$,
\be
\label{Cn-2}
T''_i\mapsto T_{i+1}\,,\quad i=1,\allowbreak\ldots,\allowbreak n-3\,,\qquad
T''_0\mapsto T_0T_1T_0T_1\,,\quad
T''_{n-2}\mapsto T_{n-1}T_nT_{n-1}T_n\,,
\ee
where $T''_0,\allowbreak\ldots,\allowbreak T''_{n-2}$ are the generators of $B(C_{n-2}^{(1)})$,
and then the quotient by the relations $T_0^2=1$, $T_n^2=1$. Then formulae
\eqref{Cn-2} induce an embedding $B(C_{n-2}^{(1)})\to\widetilde{B}(B_n^{(1)})$,
\be
\label{Bn-2}
T''_i\mapsto S_{i+1}\,,\quad i=1,\allowbreak\ldots,\allowbreak n-3\,,\qquad
T''_0\mapsto S_0S_1\,,\quad T''_{n-2}\mapsto S_{n-1}S_n\,.
\ee
Recall that \,$S_0=T_0T_1T_0$, \,$S_n=T_nT_{n-1}T_n$, and
\,$S_i=T_i$, \,$i=1,\allowbreak\ldots,\allowbreak n-1$. Formulae \eqref{Bn-2} coincide with
the embedding $\rho_0$ in (\ref{mapT2}) up to relabeling of generators.

\medskip
 \paragraph{7. One more automorphism of $B(C_n^{(1)})$.}

 ${}$

Consider the element $Z$, see~\eqref{Zn}. In addition to commutativity
\be
ZT_i=\,T_iZ\,, i=1,\allowbreak\ldots,\allowbreak n-1\,,\qquad ZX\,=\,XZ\,,
\ee
we also have
\be
ZT_1\dots T_n=\,T_0\dots T_{n-1}Z\,,\qquad
ZT_1\dots T_{n-1}=\,T_1\dots T_{n-1}Z\,,
\ee
that is, $ZT_0^{-1}X=XT_n^{-1}Z$ \,and
$ZT_0^{-1}XT_n^{-1}=T_0^{-1}XT_n^{-1}Z$. Consider an automorphism
\be
\label{phi}
\varphi:\widehat{B}(C_n^{(1)})\,\to\widehat{B}(C_n^{(1)})\,,\qquad \varphi(G)\,=\,ZGZ^{-1}\,.
\ee
Then $\varphi(T_i)\,=\,T_i$, \,$i=1,\allowbreak\ldots,\allowbreak n-1$,
\be
\label{phiT0}
\varphi(T_0)\,=\,T_0\,T_1\dots T_{n-1}\,T_n\,T_{n-1}^{-1}\dots T_0^{-1}=\,
XT_nX^{-1}.
\ee
\be
\label{phiTn}
\varphi(T_n)\,=\,T_{n-1}^{-1}\dots T_1^{-1}T_0\,T_1\dots T_{n-1}=\,
T_nX^{-1}T_0XT_n^{-1}=\,J_nT_n^{-1},
\ee
Notice that $Z$ commutes with $I_1,\allowbreak\ldots,\allowbreak I_n$ given by~\eqref{Ii}, that is,
$\varphi(I_i)=I_i$, \,$i=1,\allowbreak\ldots,\allowbreak n$.

\smallskip
The subgroup $B(C_n^{(1)})$ is invariant under the automorphism $\varphi$.


\vspace{0.2cm}

\section{$R$-matrix representation of $B_n(C^{(1)})$.\label{sec3}}
\setcounter{equation}0

Define an $R$-operator
 acting in the tensor product
 $V \otimes V$  of two $N$-dimensional vector spaces $V$
 \be
 \lb{Zam02}
 R(x,y) \cdot (\vec{e}_{k_1} \otimes \vec{e}_{k_2}) =
 (\vec{e}_{i_1} \otimes \vec{e}_{i_2}) \; R^{i_1i_2}_{k_1k_2}(x,y)  \; .
 \ee
Here vectors $\{ \vec{e}_1,\dots,\vec{e}_N \}$ form a basis in $V$ and components $R^{i_1i_2}_{k_1k_2}(x,y)$
are functions of two spectral parameters
 $x$ and $y$. Let operator $R$ satisfies Yang-Baxter equation:
 \be
 \lb{Zam04}
 R_{12}(x,y) \; R_{13}(x,z) \; R_{23}(y,z) = R_{23}(y,z) \; R_{13}(x,z) \; R_{1 2}(x,y) \;\;
 \in \;\; {\rm End}(V \otimes V \otimes V)    \; ,
 \ee
 where we have used the standard matrix notations \cite{FRT}.
 Now we introduce two matrices $||K^i_j|| \in {\rm Mat}(V)$
 and $||\overline{K}^i_j|| \in {\rm Mat}(V)$ with
 elements which are operators acting in the spaces $\widetilde{V}$
 and $\widetilde{V}'$, respectively.
 In other words we have two operators
 $K \in {\rm End}(V \otimes \widetilde{V})$ and
 $\overline{K} \in {\rm End}(V \otimes \widetilde{V}')$. Let
 these operators be solutions of
 the equation
  \be
 \lb{Zam06}
 R_{12}(x,y) \; K_1(x) \; R_{21}(y,\bar{x}) \; K_{2}(y) =
 K_{2}(y) \; R_{12}(x,\bar{y}) \; K_1(x) \; R_{21}(\bar{y},\bar{x})
 \; \in \; {\rm End}(\tilde{V}  \otimes  V \otimes V)   \; ,
 \ee
 which is called reflection equation and
 equation (cf. (\ref{Zam06}))
 \be
 \lb{Zam13}
 R_{12}(x,y) \; \overline{K}_2(y) \; R_{21}(\tilde{y},x) \; \overline{K}_{1}(x) =
 \overline{K}_{1}(x) \; R_{12}(\tilde{x},y) \; K_2(y) \; R_{21}(\tilde{y},\tilde{x})
 \; \in \; {\rm End}(V \otimes V \otimes \tilde{V}')   \; ,
 \ee
 which is called dual reflection equation.
 We explain this terminology and
 the meaning of the equations (\ref{Zam04}), (\ref{Zam06})
 and (\ref{Zam13}) in the next Section. In equations
(\ref{Zam06}) and (\ref{Zam13}) we have used notations
 \be
 \lb{invo2}
 \bar{x} = \sigma(x) \; , \;\;\;\;
 \tilde{x} = \bar{\sigma}(x) \; ,
 \ee
 where $\sigma$ and $\bar{\sigma}$ are the same
 involutive mappings $\mathbb{C} \to \mathbb{C}$
 which were introduced in (\ref{Affbg1}).

  \vspace{0.3cm}

Using operator $R(x,y)$, which is defined in (\ref{Zam02}) and (\ref{Zam04}),
  we introduce the set of $R$-operators $R_{k,k+1}(x,y)$ $(k=1,\dots,n-1)$
  which act in the space $V^{\otimes n}$
 \be
 \lb{Zam08a}
 R_{k,k+1}(x,y) = I^{\otimes (k-1)} \otimes
 R(x,y)  \;  \otimes I^{\otimes (n-k-1)}   \; .
 \ee
   For us it will be also convenient to introduce operators
 \be
 \lb{Zam08}
 \hat{R}_k(x,y) \equiv \hat{R}_{k,k+1}(x,y) = I^{\otimes (k-1)} \otimes
 P \cdot R(x)  \;  \otimes I^{\otimes (n-k-1)}  \;\;\;\;\; (k=1,\dots,n-1) \; ,
 \ee
 \be
 \lb{Zam08b}
 R_{k,r}(x,y) = P_{r,k+1} \cdot (I^{\otimes (k-1)} \otimes
 R(x)  \;  \otimes I^{\otimes (n-k-1)}) \cdot P_{r,k+1}   \; ,
 \ee
 where $P$ is a permutation operator in $V \otimes V$
 $$
 P \cdot (v_1 \otimes v_2) = (v_2 \otimes v_1) \;\;\;\;\; \forall v_1,v_2 \in V \; ,
 $$
 and $P_{r,k}= P_{k,r}$ is the permutation operator in $V^{\otimes n}$ such that
 $$
 P_{r,k} (v_1 \otimes \cdots \otimes v_k \otimes \cdots \otimes v_r \otimes \cdots \otimes v_{n}) =
 (v_1 \otimes \cdots \otimes v_r \otimes \cdots \otimes v_k \otimes \cdots \otimes v_{n}) \; .
 $$
 In terms of operators (\ref{Zam08}) equations (\ref{Zam04}), (\ref{Zam06}) and (\ref{Zam13})
 can be written in the form
 \be
 \lb{Zam09}
  \hat{R}_{k}(x,y) \; \hat{R}_{k+1}(x,z) \; \hat{R}_{k}(y,z) =
  \hat{R}_{k+1}(y,z) \; \hat{R}_{k}(x,z) \; \hat{R}_{k+1}(x,y)    \; ,
 \ee
 \be
 \lb{Zam10}
 \hat{R}_{12}(x,y) \; K_1(x) \; \hat{R}_{12}(y,\bar{x}) \; K_{1}(y) =
 K_{1}(y) \; \hat{R}_{12}(x,\bar{y}) \; K_1(x) \; \hat{R}_{12}(\bar{y},\bar{x})    \; ,
 \ee
 \be
 \lb{Zam10a}
 \hat{R}_{12}(x,y) \; \overline{K}_2(y) \; \hat{R}_{12}(\tilde{y},x) \; \overline{K}_{2}(x) =
 \overline{K}_{2}(x) \; \hat{R}_{12}(\tilde{x},y) \; \overline{K}_2(y) \; \hat{R}_{12}(\tilde{y},\tilde{x})    \; ,
 \ee
 where
 $$
 \begin{array}{c}
 K_k(x) = I^{\otimes (k-1)} \otimes K(x) \otimes I^{\otimes n-k-1} \; , \;\;\;
 \overline{K}_{k}(x) = I^{\otimes (k-1)} \otimes \overline{K}(x)
 \otimes I^{\otimes n-k-1}   \;\;\;\;\;\; 
  (k =1,\dots,n-1) \; .
 \end{array}
 $$

 Introduce the set of spectral parameters $\{ z_1 , \dots ,  z_n  \}$.
 By using the group of the elements $s_i$
 (see (\ref{Affbg1})) and matrices $\hat{R}_{k}(z_k,z_{k+1})$,
 $K_1(z_1)$, $\overline{K}_n(z_n)$ we construct the representation $\rho$ of the
 affine group $B_n(C^{(1)})$ in  $\tilde{V} \otimes V^{\otimes n} \otimes  \tilde{V}'$
 \be
 \lb{Zam10b}
 \begin{array}{c}
 \rho(T_i) \; = \;  s_i \, \hat{R}_{i}(z_i,z_{i+1})  \;\;\;\; (i=1,\dots,n-1) \; , \;\;\;
 \rho(T_0) \; = \;   K_1(z_1) \, s_0  \; , \;\;\;
 \rho(T_n) \; = \;  \overline{K}_n(z_n) \, s_n  \; .
 \end{array}
 \ee
 One can  directly check that $\rho(T_i)$ $(i=0,\dots,n)$ satisfy
 defining relations in (\ref{Affbg}), (\ref{Affbg2}) if $\hat{R}_{k}(z_k,z_{k+1})$
 and  $K_1(z_1)$, $\overline{K}_n(z_n)$
  satisfy relations (\ref{Zam09}), (\ref{Zam10}), (\ref{Zam10a}).

   Further we will use the operator
 $D_{z_k}$ such that for any wave function $\Psi(z_1, \dots , z_n)$
 and any operator $f(z_1, \dots , z_n)$ we have
 \be
\lb{Zam57}
\begin{array}{c}
 D_{z_k} \cdot f(z_1, \dots , z_k, \dots , z_n) =
 f(z_1, \dots , \tilde{\bar{z_k}}, \dots , z_n) \cdot D_{z_k} \; , \\ [0.2cm]
 D_{z_k} \cdot \Psi(z_1, \dots , z_k, \dots , z_n) =
 \Psi(z_1, \dots , \tilde{\bar{z_k}}, \dots , z_n)  \; ,
 \end{array}
 \ee
 where $\tilde{\bar{z_k}} = \bar{\sigma} (\sigma (z_k))$.
 We note that the operator $D_{z_k}$
 in (\ref{Zam57}) can be written in the representation (\ref{Affbg1}) as
 \be
\lb{Zam57D}
 D_{z_k} = (s_{k-1} \cdots s_1) (s_0 \cdots s_n) (s_{n-1} \cdots s_k) = {\sf s}(J_k) \; ,
 \ee
 where elements $J_k$ were introduced in (\ref{jucys1}).

\vspace{0.2cm}

  \noindent
  {\bf Theorem \ref{sec3}.1.} {\it The images
  of the commutative elements (\ref{jucys1}) are
  operators in $\tilde{V} \otimes V^{\otimes n} \otimes \tilde{V}'$
  \be
  \lb{Zam05ww}
   \begin{array}{c}
  \rho(J_i)={\sf A}_i =  \hat{R}_{k-1}^{-1}(z_{k-1},z_k) \cdots  \hat{R}_1^{-1}(z_1,z_k) \; K_1(z_k) \; \hat{R}_{1}(z_1,\bar{z}_k) \cdots \hat{R}_{k-1}(z_{k-1},\bar{z}_k) \cdot \\ [0.2cm]
  \hat{R}_{k}(z_{k+1},\bar{z}_k) \cdots \hat{R}_{n-1}(z_{n},\bar{z}_k) \;
  \overline{K}_n(\bar{z}_k) \;
 \cdot D_{z_k} \; \cdot
\hat{R}_{n-1}(z_k,z_n)  \cdots \hat{R}_{k}(z_k,z_{k+1}) \; ,
  \\ [0.3cm]
  \rho(\overline{J}_i)=\overline{\sf A}_i = \hat{R}_{k-1}(z_k,z_{k-1}) \cdots  \hat{R}_1(z_k,z_1) \; K_1(z_k) \; \hat{R}_{1}(z_1,\bar{z}_k) \cdots \hat{R}_{k-1}(z_{k-1},\bar{z}_k) \cdot \\ [0.2cm]
  \hat{R}_{k}(z_{k+1},\bar{z}_k) \cdots \hat{R}_{n-1}(z_{n},\bar{z}_k) \;
  \overline{K}_n(\bar{z}_k) \;
 \cdot D_{z_k} \; \cdot
\hat{R}_{n-1}^{-1}(z_n,z_k)  \cdots
\hat{R}_{k}^{-1}(z_{k+1},z_k) \; ,
  \end{array}
  \ee
  form two sets of flat connections for
  quantum Knizhnik-Zamolodchikov
  equations
  \be
 \lb{Zam16w}
 \begin{array}{c}
 {\sf A}_k(z_1, \dots , z_k, \dots , z_n) \;
   \; \Psi(z_1, \dots , z_k, \dots , z_n) = \Psi(z_1, \dots , z_k, \dots , z_n) \; , \\ [0.2cm]
   \overline{\sf A}_k(z_1, \dots , z_k, \dots , z_n) \;
   \; \overline{\Psi}(z_1, \dots , z_k, \dots , z_n) =
   \overline{\Psi}(z_1, \dots , z_k, \dots , z_n) \; ,
   \end{array}
 \ee
 where functions
 $\Psi, \overline{\Psi} \in
 \tilde{V} \otimes V^{\otimes n} \otimes \tilde{V}'$. } \\
 {\bf Proof.} Formulas (\ref{Zam05ww}) are obtained by direct calculations. The flatness of the connections (\ref{Zam16w})
 $$
 [{\sf A}_k , \, {\sf A}_j] = 0 =
 [ \overline{\sf A}_k , \, \overline{\sf A}_j] \; ,
 $$
 follows from the Proposition \ref{ABG}.1. \hfill \qed

 \section{Flat connections for quantum Knizhnik-Zamolodchikov equations. Approach with Zamolodchikov algebra.\label{sec4}}
\setcounter{equation}0

\subsection{Zamolodchikov algebra.}

Introduce a set of operators $A^i(z)$ $(i=1,2,\dots,N)$ which act in the
 complex vector space ${\cal H}$.
Each operator $A^i(z)$ is a function of the spectral parameter $z$.
The operators $A^i(z)$ are generators of the algebra ${\cal Z}$
 with quadratic defining relations
 (see e.g.  \cite{GoshZam} and references therein)
 \be
 \lb{Zam01}
 A^{i_1}(x) \; A^{i_2}(y) = R^{i_1i_2}_{k_1k_2}(x,y) \; A^{k_2}(y) \; A^{k_1}(x)   \; .
 \ee
 where $R^{i_1i_2}_{k_1k_2}(x,y) \in \mathbb{C}$ are functions
 of the spectral parameters $x$ and $y$ and also
 are components of an $R$-operator
 acting in the tensor product
 $V \otimes V$  of two $N$-dimensional vector spaces $V$
 (see (\ref{Zam02})).
The algebra ${\cal Z}$ is called Zamolodchikov algebra.
Relations (\ref{Zam01}) can be written in concise matrix notations \cite{FRT} as following
\be
 \lb{Zam03}
 A^{1 \rangle}(x) \; A^{2 \rangle}(y) = R_{1 2}(x,y) \; A^{2 \rangle}(y) \; A^{1 \rangle}(x)   \; .
 \ee
 Consider the product $A^{i_1}(x) A^{i_2}(y) A^{i_3}(z)$ of three operators and reorder it
  with the help of (\ref{Zam01}) as following
 $$
 A^{i_1}(x) A^{i_2}(y) A^{i_3}(z)  \;\; \to \;\; A^{k_3}(z) A^{k_2}(y) A^{k_1}(x) \; ,
 $$
 in two different ways in accordance with the arrangement of brackets
 \be
 \lb{lexy}
 (A^{i_1}(x) A^{i_2}(y)) A^{i_3}(z) = A^{i_1}(x) (A^{i_2}(y) A^{i_3}(z)) \; .
 \ee
  As a result we obtain the
 self-consistence condition for the matrix $R(x,y)$
 in the form of the Yang-Baxter equation (\ref{Zam04}).
The solutions $R(x,y)$ of the equation (\ref{Zam04})
 define Zamolodchikov algebra (\ref{Zam01}).

 Now we extend (see \cite{GoshZam})
 the algebra ${\cal Z}$ by adding new ''boundary" operators $B^\alpha$
 ($\alpha =1,2, \dots,M$) which act in ${\cal H}$ and obey  relations
 \be
 \lb{Zam05}
 \begin{array}{c}
 A^{i}(x) \; B^\alpha = K^{i \alpha}_{k \beta}(x) \; A^{k}( \bar{x} ) \; B^\beta
  \;\;\; \Rightarrow \;\;\;
 A^{1 \rangle}(x) \; B = K_{1}(x) \; A^{1 \rangle}( \bar{x} ) \; B  \; , \\ [0.2cm]
 \bar{x} = \sigma(x) \in \mathbb{C} \; ,
 \end{array}
 \ee
 where $\bar{x}$ is a reflected spectral parameter
 and $\sigma$ -- involutive operation $\mathbb{C} \to \mathbb{C}$ such that $\sigma^2 =1$. E.g., for
  rational and trigonometric $R$-matrices (pay attention to the special dependence
  of spectral parameters)
  \be
  \lb{RatTri}
  R(x,y)=R(x-y) \; , \;\;\; R(x,y)=R(x/y)  \; ,
  \ee
  one can take
  $\sigma = \sigma_a$ and $\sigma = \sigma^{tri}_b$, respectively, where
   \be
 \lb{Zam05s}
   \sigma_a (x) = a - x \; , \;\;\;\; \sigma^{tri}_b(x) = b/x \; ,
   \ee
   and $a,b \in \mathbb{C}$ are parameters which specify involutions $\sigma$, $\sigma^{tri}$.
   Matrix $K$ with components $K^{i \alpha}_{k \beta}(x)$ acts in the space $V \otimes \widetilde{V}$,
 where $\widetilde{V}$ is $M$-dimensional vector space. This matrix
 is called reflection matrix and describes a reflection of particles
 from right boundary \cite{GoshZam}. For simplicity, in the second formula in (\ref{Zam05}) and
 below, we omit indices $\alpha,\beta,\dots$ related to the space $\widetilde{V}$.

 In the same way as in (\ref{lexy}), one can
 consider two different ways for the reordering of special product of
 3 generators (including $B^\alpha$):
 $$
 [A^{i_1}(x) \; A^{i_2}(y)] \; B =  A^{i_1}(x) \; [A^{i_2}(y) \; B ]
 \;\;\;\; \longrightarrow \;\;\;\;
  A^{k_1}(\bar{x}) \; A^{k_2}(\bar{y}) \; B  \; .
 $$
 As a result, in addition to the Yang-Baxter equation (\ref{Zam04}),
  we obtain new consistence condition for the reflection matrix $K$ in the form of the reflection equation (\ref{Zam06}).

 Now besides the "right" boundary operators $B^\alpha$ with relations
 (\ref{Zam05}), we  also introduce the  "left" boundary operators $\overline{B}_{\alpha'}$
 $(\alpha' = 1,\dots,M')$ with relations
 \be
 \lb{Zam12}
 \begin{array}{c}
 \overline{B}_{\beta'} \; A^{i}(x) =  \overline{K}^{i\alpha'}_{k \beta'}(x) \; \overline{B}_{\alpha'}\;
 A^{k}( \tilde{x} )   \;\;\; \Rightarrow \;\;\;
     \overline{B} \; A^{1 \rangle}(x)  =
      \overline{K}_{1}(x) \; \overline{B} \; A^{1 \rangle}( \tilde{x} )  \; , \\ [0.2cm]
      \tilde{x} = \bar{\sigma}(x)  \in \mathbb{C} \; ,
      \end{array}
 \ee
 where $\bar{\sigma}$ is another
 involutive operation in $\mathbb{C}$: $\bar{\sigma}^2=1$
 (e.g., one can define $\bar{\sigma}$ as in (\ref{Zam05s}) but with another
 parameters $a,b$). Equations (\ref{Zam12}) and operator
 $\overline{K}(x) \in {\rm End}(V \otimes \tilde{V}')$,
  where $\widetilde{V}'$ is $M'$-dimensional vector space, describe the reflection
 of particles from the left boundary. Two different ways for the reordering of the
  product of 3 generators (including $\overline{B}$):
 $$
 \overline{B} \; A^{i_1}(x) \; A^{i_2}(y)   \;\;\;\; \longrightarrow \;\;\;\;
   \overline{B} \; A^{k_1}(\tilde{x}) \; A^{k_2}(\tilde{y})   \; ,
 $$
 give additional consistence condition
 in the form of the dual reflection equation (\ref{Zam13}).

 Note that applying  defining relations (\ref{Zam03}),
 (\ref{Zam05}) and (\ref{Zam12}) twice, we deduce three unitary relations for matrices $R$, $\overline{K}$ and $K$
  \be
 \lb{Zam07}
 R_{12}(x,y) \; R_{21}(y,x)  = I \otimes I \; , \;\;\;
 K_{1}(x)  \; K_1(\bar{x})  = I \otimes \widetilde{I}  \; , \;\;\;
 \overline{K}_{1}(x)  \;  \overline{K}_1(\tilde{x})  = I \otimes \widetilde{I}'  \; ,
 \ee
 where $I$, $\widetilde{I}$  and $\widetilde{I}'$
 -- unite operators in $V$, $\widetilde{V}$ and $\widetilde{V}'$, correspondingly.

  In physics the matrices $R$ and $K$, $\overline{K}$ which satisfy equations
   (\ref{Zam04}), (\ref{Zam06}), (\ref{Zam13}) and
 (\ref{Zam07}) describe the factorizable scattering on a half line \cite{GoshZam}, \cite{Cher1}, or
 define the integrable spin chains with nontrivial boundary conditions \cite{Skl}.

  Note that
 if matrices $R$, $K$ and $\overline{K}$ satisfy
 unitarity conditions (\ref{Zam07}), then
 for the representation (\ref{Zam10b}) we have
 $(\rho(T_i))^2 = I$, where $I$ is the unit operator in
 $\tilde{V} \otimes V^{\otimes n} \otimes  \tilde{V}'$. Thus, in this case the equations (\ref{Zam10b})
 define the representation of the Coxeter group $W_n(C^{(1)})$.

 \subsection{Flat connections for quantum Knizhnik-Zamolodchikov equations.}

 Consider the boundary Zamolodchikov algebra ${\cal Z}_{LR}$ with generators
 $\{ A^i(x), \; B^\alpha, \; \overline{B}_{\beta'} \}$.
 Namely, the algebra ${\cal Z}_{LR}$ includes the generators
 $A^i(x)$ of the Zamolodchikov algebra ${\cal Z}$ and both left and right boundary
 operators $B^\alpha$ and  $\overline{B}_{\beta'}$.
  Consider the special element in ${\cal Z}_{LR}$:
 \be
 \lb{Zam14}
 [\Psi^\alpha_{\beta'}]^{i_n \dots i_1}(z_n, \dots , z_k, \dots , z_1) =
 \overline{B}_{\beta'}  \; A^{i_n}(z_n) \cdots A^{i_k}(z_k)
 \cdots  A^{i_1}(z_1)  \; B^\alpha \; ,
 \ee
 and push the $k$-th operator $A^{i_k}(z_k)$, in the ordered product
 $\bigl( A^{i_n}(z_n) \cdots  A^{i_1}(z_1)\bigr)$ in the right hand side of (\ref{Zam14}),
 with the help of equations (\ref{Zam01}) to the right. Then we reflect this
 operator  from the right boundary
 operator $B^\alpha$ with the help of (\ref{Zam05}), and push
 the reflected operator $A_{(k)}(\bar{z}_k)$ backward to the left
 with the help of (\ref{Zam01}) up to
 the left boundary operator $\overline{B}_{\beta'}$. Then we reflect the
  operator $A_{(k)}(\bar{z}_k)$ from this
 boundary operator and finally place the operator $A_{(k)}(\tilde{\bar{z}}_k)$
  on its initial $k$-th position in the ordered
 product $A_{(n)}(z_n) \cdots A_{(2)}(z_2) \;  A_{(1)}(z_1) $.
 As a result we obtain the equation
  \be
 \lb{Zam15}
 \begin{array}{c}
 (\Psi^\alpha_{\beta'})^{i_1 \dots i_n}(z_n, \dots , z_k, \dots , z_1) 
 = [{\cal A}_k(z_1, \dots , z_k, \dots , z_n)]^{i_1 \dots i_n; \; \alpha \gamma'}_{j_1 \dots j_n;
  \; \beta' \delta}
 \;\; (\Psi^\delta_{\gamma'})^{j_1 \dots j_n}
 (z_n, \dots , \tilde{\bar{z}}_k, \dots , z_1) \; ,  \\ [0.2cm]
 \tilde{\bar{z}}_k=\bar{\sigma} (\sigma (z_k)) \; ,
 \end{array}
 \ee
 where involutions $\sigma$ and $\bar{\sigma}$ were introduced in
 (\ref{Zam05}) and (\ref{Zam12}) while the matrix
 $$
 \begin{array}{c}
 [{\cal A}_k(z_1, \dots , z_k, \dots , z_n)]_{12 \dots n} =
 K_{k}(z_k; \, \vec{z}_{(1,k-1)}) \cdot \overline{K}_{k}(\bar{z}_k; \, \vec{z}_{(k+1,n)}) \; , \\ [0.3cm]
 \vec{z}_{(1,k-1)}=(z_1,\dots,z_{k-1}) \; ,
 \;\;\; \vec{z}_{(k+1,n)}=(z_{k+1},\dots,z_{n}) \; ,
 \end{array}
 $$
 is defined by means of dressed reflection matrices
\be
\lb{Zam55}
\begin{array}{c}
K_{k}(x;\vec{z}_{(1,k-1)}) =
R_{k,k-1}(x,z_{k-1}) \cdots  R_{k1}(x,z_1) \; K_k(x) \;
R_{k1}(z_1,\bar{x})  \cdots R_{k,k-1}(z_{k-1},\bar{x}) = \\ [0.2cm]
 = \hat{R}_{k-1}^{-1}(z_{k-1},x) \cdots \hat{R}_2^{-1}(z_2,x) \, \hat{R}_1^{-1}(z_1,x) \; K_1(x) \;
\hat{R}_{1}(z_1,\bar{x}) \, \hat{R}_2(z_2,\bar{x}) \cdots \hat{R}_{k-1}(z_{k-1},\bar{x}) = \\ [0.2cm]
= \hat{R}_{k-1}^{-1}(z_{k-1},x) \, K_{k-1}(x;\vec{z}_{(k-2)}) \, \hat{R}_{k-1}(z_{k-1},\bar{x}) \; ,
\end{array}
\ee
\be
\lb{Zam56}
\begin{array}{c}
\overline{K}_{k}(\bar{x};\vec{z}_{(k+1,n)}) =
R_{k+1,k}(z_{k+1},\bar{x}) \cdots R_{nk}(z_{n},\bar{x}) \; \overline{K}_k(\bar{x}) \;
 R_{kn}(\tilde{\bar{x}},z_n)  \cdots R_{k,k+1}(\tilde{\bar{x}},z_{k+1}) = \\ [0.2cm]
\hat{R}_{k}(z_{k+1},\bar{x}) \cdots \hat{R}_{n-1}(z_{n},\bar{x}) \; \overline{K}_n(\bar{x}) \;
\hat{R}_{n-1}(\tilde{\bar{x}},z_n)  \cdots \hat{R}_{k}(\tilde{\bar{x}},z_{k+1}) = \\ [0.2cm]
= \hat{R}_{k}(z_{k+1},\bar{x}) \,
\overline{K}_{k+1}(\bar{x};\vec{z}_{(k+2,n)}) \, \hat{R}_{k}(\tilde{\bar{x}},z_{k+1}) \; .
\end{array}
\ee
To write expression (\ref{Zam55}) for the matrix $K_{k}(x;\vec{z}_{(1,k-1)})$ we
take into account the unitarity condition for the $R$-operator
 $\hat{R}_{k}(x,z) =\hat{R}_{k}^{-1}(z,x)$.

 For rational and trigonometric $R$-matrices
 (\ref{RatTri}) the involutions $\sigma$ and $\bar{\sigma}$ could be defined as in
 (\ref{Zam05s})
 $$
 {\rm rational \;\; case}: \;\; \sigma = \sigma_a \; , \;\;\;  \bar{\sigma} = \sigma_{a'} \;\; ; \;\;\;
  {\rm trigonometric \;\; case}: \;\; \sigma = \sigma^{tri}_b \; , \;\;\;
    \bar{\sigma} = \sigma^{tri}_{b'} \; ; \;\;\;
 $$
 and  we respectively obtain
 \be
 \lb{Zam05ss}
 \tilde{\bar{x}} = \sigma_{a'}( \sigma_a ( x)) = (a'-a) + x \; , \;\;\;\;
 \tilde{\bar{x}} = \sigma^{tri}_{b'} ( \sigma^{tri}_b ( x)) = \frac{b'}{b} \, x \; ,
 \ee
 i.e., for the rational case the spectral
 parameter $\tilde{\bar{x}}$ is a shift of $x$ by a constant $(a'-a)$, while
 for the trigonometric case the parameter $\tilde{\bar{x}}$ is a multiplication of
 $x$ by a constant $b'/b$. In view of this, for rational and trigonometric cases the operator
 $D_z$ (\ref{Zam57}), (\ref{Zam57D}) can be considered as finite difference derivatives. Note that
 $\bar{\sigma} \sigma \neq \sigma \bar{\sigma} $.

 One can write eqs. (\ref{Zam15}) in the form of quantum Knizhnik-Zamolodchikov
 equations (see (\ref{Zam16w}):
 \be
 \lb{Zam16}
 {\sf A}_k(z_1, \dots , z_k, \dots , z_n) \;
   \; \Psi(z_1, \dots , z_k, \dots , z_n) = \Psi(z_1, \dots , z_k, \dots , z_n) \; ,
 \ee
 where we interpret $\Psi$ (\ref{Zam14}) as a wave function and introduce connections
  \be
 \lb{Zam05f}
 \begin{array}{c}
 {\sf A}_k(z_1, \dots , z_k, \dots , z_n) = {\cal A}_k(z_1, \dots , z_k, \dots , z_n) \; D_{z_k} =
 K_{k}(z_k; \, \vec{z}_{(1,k-1)}) \cdot \overline{K}_{k}(\bar{z}_k; \, \vec{z}_{(k+1,n)})
 \; D_{z_k} = \\ [0.3cm]
 = K_{k}(z_k; \, \vec{z}_{(1,k-1)}) \cdot \overline{\sf K}_{k}(\bar{z}_k; \, \vec{z}_{(k+1,n)}) \; .
  \end{array}
 \ee
In the right hand side of (\ref{Zam05f}) we use the dressed reflection matrix (\ref{Zam55})
for $x=z_k$ which can be written in the representations
(\ref{Affbg1}) and (\ref{Zam10b}) as the following
 \be
 \lb{Zam05d}
  \begin{array}{c}
 K_{k}(z_k; \, \vec{z}_{(1,k-1)})  = \hat{R}_{k-1}^{-1}(z_{k-1},z_k) \cdots  \hat{R}_1^{-1}(z_1,z_k) \; K_1(z_k) \; \hat{R}_{1}(z_1,\bar{z}_k) \cdots \hat{R}_{k-1}(z_{k-1},\bar{z}_k) = \\ [0.2cm]
=  \rho(T_{k-1}^{-1} \cdots T_1^{-1} \, T_0 \, T_1 \cdots T_{k-1}) \cdot
 (s_{k-1} \cdots s_1 \, s_0 \, s_1 \cdots s_{k-1}) = \rho(\bar{a}_k) \cdot {\bf s}(a_k) \; ,
  \end{array}
 \ee
 where $\bar{a}_k = T_{k-1}^{-1} \cdots T_1^{-1} T_0 T_1 \cdots T_{k-1}$
 and elements $a_k$ were defined in (\ref{jucys1}). Besides this we also
 define new dressed reflection matrix
  \be
 \lb{Zam05K}
 \begin{array}{c}
 \overline{\sf K}_{k}(\bar{x}; \, \vec{z}_{(k+1,n)}) =
 \overline{K}_{k}(\bar{x}; \, \vec{z}_{(k+1,n)})
 \; D_{x} = \hat{R}_{k}(z_{k+1},\bar{x}) \cdot
   \overline{\sf K}_{k+1}(\bar{x}; \, \vec{z}_{(k+2,n)}) \cdot \hat{R}_{k}(x,z_{k+1}) =
   \\ [0.3cm]
 = \hat{R}_{k}(z_{k+1},\bar{x}) \cdots \hat{R}_{n-1}(z_{n},\bar{x}) \; \overline{K}_n(\bar{x}) \;
 \cdot D_{x} \; \cdot
\hat{R}_{n-1}(x,z_n)  \cdots \hat{R}_{k}(x,z_{k+1}) \; ,
\end{array}
 \ee
 which includes the finite difference operator $D_{x}$ (\ref{Zam57}).
  In the representations
(\ref{Affbg1}) and (\ref{Zam10b}), for $x=z_k$, the
 matrix (\ref{Zam05K}) can be written  as the following
\be
 \lb{Zam05dd}
 \overline{\sf K}_{k}(\bar{z}_k; \, \vec{z}_{(k+1,n)}) =
 (s_{k-1} \cdots s_1 \, s_0 \, s_1 \cdots s_{k-1}) \cdot \rho( T_{k} \cdots
  T_{n-1} \, T_n \, T_{n-1} \cdots  T_{k}) = {\bf s}(a_k) \cdot \rho(b_k)  \; ,
 \ee
 where $a_k$ and $b_k$ were defined in (\ref{jucys1}).
To obtain relations (\ref{Zam05d}) and (\ref{Zam05dd}) we have used formulas
(\ref{Zam57D}) and
 $$
 \begin{array}{c}
  \bar{z}_k = (s_{k-1} \cdots s_1 \, s_0 \, s_1 \cdots s_{k-1}) \, z_k \,
  (s_{k-1} \cdots s_1 \, s_0 \, s_1 \cdots s_{k-1}) \; .
  \end{array}
  $$
  Finally, using (\ref{Zam05d}) and (\ref{Zam05dd})
  one can write connections (\ref{Zam05f}) in the form
   \be
 \lb{Zam05ff}
   {\sf A}_k(z_1,
   \dots , z_n) =  \rho(\bar{a}_k) \cdot \rho(b_k) = \rho(J_k) \; .
  \ee

 Applying equation (\ref{Zam16}) twice (for two different indices $k$ and $r$)
 we deduce the consistency condition
 $$
 [{\sf A}_k , \; {\sf A}_r] \;\; \Psi(z_1, \dots , z_n) = 0 \; ,
 $$
 and our conjecture is that the connections $A_k$,
 explicitly given in (\ref{Zam05f}) and (\ref{Zam05ff}), are flat:
\be
\lb{Zam17}
 [{\sf A}_k , \; {\sf A}_r] = 0  \; .
 \ee
 One can prove this identity directly by using the fact
 that connections $A_k$ (\ref{Zam05ff}) are the images of the commuting elements
 $J_k \in B_n(C^{(1)})$ (see Proposition \ref{ABG}.1). Note that commutativity (\ref{Zam17}) of
 connections $A_k$ (\ref{Zam05f}), where matrix $K_{k}(z_k; \, \vec{z}_{(1,k-1)})$ is taken
 in the form (\ref{Zam05d}), is valid even
 for the case when $R$-matrix is not satisfies unitarity condition. So, we have proved
 the following statement:

   \vspace{0.2cm}

 \noindent
 {\bf Theorem \ref{sec4}.1.} {\it Connections ${\sf A}_k$ which were defined in (\ref{Zam05f}),
 (\ref{Zam05d}), (\ref{Zam05K}) are flat (\ref{Zam17})
  for any matrices $R$, $K$ and $\overline{K}$ satisfying eqs.
  (\ref{Zam09}),(\ref{Zam10}) and (\ref{Zam10a}) and
  any involutive operations $\sigma,\bar{\sigma}$.}

\vspace{0.2cm}

 \noindent
{\bf Remark 1.} One can think about boundary operators $B^\alpha$ and $\overline{B}_{\alpha'}$
in (\ref{Zam05}),
(\ref{Zam12}) and (\ref{Zam14}) as about boundary states $| B^\alpha \rangle \in {\cal H}$ and
$\langle \overline{B}_{\alpha'} | \in {\cal H}^*$
with the same conditions as in (\ref{Zam05}), (\ref{Zam12}).
In this case the operator (\ref{Zam14}) is represented as the matrix element
 \be
 \lb{Zam18}
[\Psi^\alpha_{\beta'}]^{i_n \dots i_1}(z_n, \dots , z_2,  z_1) =
 \langle \overline{B}_{\beta'}  | \; A^{i_n}(z_n) \cdots A^{i_2}(z_2) \;
  A^{i_1}(z_1)  \; | B^\alpha \rangle \; ,
 \ee
and the equation (\ref{Zam16}), with the wave function $\Psi$ which
is given in (\ref{Zam18}), is nothing but the
quantum Knizhnik-Zamolodchikov (q-KZ) equations for the
system with nontrivial boundary conditions.  One can put
$\tilde{V} = \tilde{V}'$, $\beta' = \alpha$ in (\ref{Zam18}) and sum over $\alpha$.
As a result we obtain the following form of the solution of q-KZ equation
 \be
 \lb{Zam19}
 \Psi^{i_n \dots i_1}(z_n, \dots , z_2,  z_1) = {\rm Tr}_{\cal H}
 \Bigl( \; A^{i_n}(z_n) \cdots A^{i_2}(z_2) \;
  A^{i_1}(z_1)  \; \rho \, \Bigr) \; ,
 \ee
 where $\rho = | B^\alpha \rangle \langle \overline{B}_\alpha  |$ can be
 considered as a density matrix.

\vspace{0.2cm}

 \noindent
{\bf Remark 2.} For systems with periodic boundary conditions one can deduce q-KZ equations
by using the same method as was used above for the systems with
nontrivial boundary conditions and open boundaries. Consider
the function (\ref{Zam19}) with any operator  $\rho$ and
 require that this operator satisfies
 commutation relations with generators $A^i(x)$:
 \be
 \lb{Zam20}
 A^i(x) \,  \rho = Q^i_j(x) \, \rho \, A^i(\tilde{\bar{x}})  \; , \;\;\;\;
 \tilde{\bar{x}} = \bar{\sigma} (\sigma (x)) \; .
 \ee
Here functions $Q^i_j(x)$ are components of a numerical matrix. Taking into account (\ref{Zam20}) we obtain the following periodicity condition
 for the wave function (\ref{Zam19})
 \be
 \lb{percond}
 \begin{array}{c}
 \Psi^{i_n \dots i_1}(z_n, \dots , z_2,  z_1) =  \\ [0.2cm]
  = {\rm Tr}_{\cal H}
 \Bigl( \; A^{i_n}(z_n) \cdots A^{i_3}(z_3) \;
  A^{i_2}(z_2) \, Q^{i_1}_{j_1}(z_1) \; \rho \; A^{j_1}(\tilde{\bar{z}}_1)  \, \Bigr) =
  Q^{i_1}_{j_1}(z_1) \, \Psi^{j_1 i_n \dots i_2}( \tilde{\bar{z}}_1, z_n, \dots , z_2) \; .
  \end{array}
 \ee
 The associativity equation
  $A^{i_1}(x)(A^{i_2}(y)\; \rho) = (A^{i_1}(x)A^{i_2}(y)\;) \rho$
 requires consistency condition for matrix $Q^i_j(x)$
  \be
 \lb{Zam21}
 R_{12}(z_1,z_2) Q_1(z_1) Q_2(z_2)  =
 Q_1(z_1) \, Q_2(z_2) \, R_{12}(\tilde{\bar{z}}_1,\tilde{\bar{z}}_2)  \; .
 \ee
 We also require the condition
 $$
 R_{12}(\tilde{\bar{z}}_1,\tilde{\bar{z}}_2) = R_{12}(z_1,z_2) \;\;\;\; \Leftrightarrow \;\;\;\;
 D_{z_1} \, D_{z_2} \, R_{12}(z_1,z_2) = R_{12}(z_1,z_2) \, D_{z_1} \, D_{z_2} \; ,
 $$
 which is obtained automatically for the rational and trigonometric cases,
 when involutions $\bar{\sigma}$, $\sigma$ are fixed as in (\ref{Zam05ss}).
 In this case equation (\ref{Zam21}) is written as
 $$
 (D_{z_1} \, Q_1(z_1)) \, (D_{z_2} \, Q_2(z_2)) \, R_{12}(z_1,z_2) =
 R_{12}(z_1,z_2) \, (Q_1(z_1) \,  D_{z_1}) \, (Q_2(z_2) \, D_{z_2}) \; .
 $$
 Now we again pick up  the generator $A^{i_k}(z_k)$  in the right hand side of
 (\ref{Zam19}) push this generator to the right with the help of (\ref{Zam01}), then
 use relation (\ref{Zam20}) and cyclic property of the trace and finally place the
 operator $A^{i_k}(\tilde{\bar{z}}_k)$ on its initial $k$-th position. As a result we obtain
 equation
 \be
 \lb{Zam22}
 \Psi(z_n, \dots , z_2,  z_1)  = {\sf A}_k(\vec{z}_{(1,n)}) \; \Psi(z_n, \dots , z_2,  z_1) \; ,
 \ee
 where $\vec{z}_{(1,n)}=(z_1,\dots,z_n)$ and ${\sf A}_k(\vec{z}_{(1,n)})$ is the flat connection for q-KZ
 equation in the periodic case \cite{FreRe}:
 \be
 \lb{Zam22f}
 \begin{array}{cl}
  {\sf A}_k(\vec{z}_{(1,n)}) =  &
  R_{k,k-1}(z_k,z_{k-1}) \cdots R_{k,2}(z_k,z_{2}) R_{k,1}(z_k,z_{1}) \,
  Q_k(z_k) \, D_{z_k} \cdot \\ [0.3cm]
  \, & \cdot R_{kn}(z_k,z_n) R_{k,n-1}(z_k,z_{n-1}) \cdots R_{k,k+1}(z_k,z_{k+1}) \; .
 \end{array}
 \ee
 Here the finite difference operator $D_{z_k}$ is the same as in (\ref{Zam57}).
 Using for the periodic
  braid group elements $T_i$ the same $R$-matrix
  representation (\ref{Zam10b}) we write connection (\ref{Zam22f}) as
  (cf. (\ref{perBn3}))
  \be
 \lb{Zam22Ap}
 \begin{array}{c}
  {\sf A}_k(\vec{z}_{(1,n)}) =
  \rho(T_{k-1} \cdots T_1) \cdot X
 \cdot \rho( T^{-1}_{n-1}  \cdots T_k^{-1}) \; , \\ [0.3cm]
  X : =  Q_1(z_1) \, D_{z_1} \, \hat{\sf s}_1 \cdots \hat{\sf s}_{n-1} \; ,
 \end{array}
 \ee
 where $\hat{\sf s}_k = P_{k,k+1} \, {\sf s}_k$ and we have used unitarity
 conditions $T_i^2=1$. We have (for simplicity  we write $T_i$ instead of
 $\rho(T_i)$)
 \be
 \lb{ZamTX2}
 \begin{array}{c}
 T_i \, \hat{\sf s}_{i+1} \, \hat{\sf s}_i =
 \hat{\sf s}_{i+1} \, \hat{\sf s}_i \, T_{i+1} \; , \;\;\;
  T_{i+1} \, \hat{\sf s}_i \, \hat{\sf s}_{i+1} =
 \hat{\sf s}_i \, \hat{\sf s}_{i+1} \, T_i \; , \\ [0.2cm]
 X \, T_i = T_{i+1} \, X \; , \;\;\; (i=1,\dots , n-2) \; , \\ [0.2cm]
 T_1 \cdot X^2 = T_1 Q_1 \, D_{z_1} Q_2 \, D_{z_2}
 (\hat{\sf s}_1 \cdots \hat{\sf s}_{n-1})^2 =
Q_1 \, D_{z_1} Q_2 \, D_{z_2} \,  T_1 \,
 (\hat{\sf s}_1 \cdots \hat{\sf s}_{n-1})^2 = \\ [0.2cm]
  =  Q_1 \, D_{z_1} Q_2 \, D_{z_2} \,  T_1 \,
 (\hat{\sf s}_2 \hat{\sf s}_1) \cdot (\hat{\sf s}_3 \hat{\sf s}_2) \cdots (\hat{\sf s}_{n-1} \hat{\sf s}_{n-2}) = \\ [0.2cm]
 = Q_1 \, D_{z_1} Q_2 \, D_{z_2} \,
 (\hat{\sf s}_2 \hat{\sf s}_1)
  \cdots (\hat{\sf s}_{n-1} \hat{\sf s}_{n-2})  \,  T_{n-1}
  = X^2 \,  T_{n-1} \; .
 \end{array}
 \ee
 One can check that the element
 \be
 \lb{Xn}
 T_n :=
 X^{-1} \, T_1 \cdot X = X \,  T_{n-1} \, X^{-1} \; ,
 \ee
  satisfies
 periodic braid relations
 $$
 T_n \,  T_{n-1} \,  T_n = T_{n-1} \, T_n \, T_{n-1} \; , \;\;\;
 T_n \,  T_{1} \,  T_n = T_{1} \, T_n \, T_{1} \; .
 $$
 Let $T_1$ be unitary operator $T_1^2 =1$. In this case
 the connection (\ref{Zam22Ap}) satisfies the periodicity condition
 $$
  {\sf A}_k(\vec{z}_{(1,n)}) =
  T_{k-1} \cdots T_1 \cdot X
 \cdot  T^{-1}_{n-1}  \cdots T_k^{-1} =
  T_{k-1} \cdots T_2 \cdot X
 \cdot  T_{n}^{-1} T^{-1}_{n-1}  \cdots T_k^{-1} \; .
 $$

\vspace{0.2cm}

 \noindent
 {\bf Proposition~ \ref{sec4}.2}~\cite{FreRe}.~ {\it
  For the periodic chain the connections (\ref{Zam22Ap})are flat,
  i.e. satisfy (\ref{Zam17}).} \\
 {\bf Proof.} The proof is the same as the proof of the Proposition \ref{ABG}.1
 in Section {\bf \ref{ABG}}.
 \hfill \qed

 \vspace{0.3cm}

 \noindent
{\bf Remark 3.}
 Consider operator $T_{V{\cal V}}(x) \in {\rm End}(V \otimes {\cal V})$ which satisfies
the intertwining relations
\be
 \lb{Zam23}
 \begin{array}{c}
 {\cal R}_{{\cal V}{\cal V}'}(x,y) \;  T_{1 {\cal V}}(x) \;  T_{1 {\cal V}'}(y) =
 T_{1 {\cal V}'}(y)\;  T_{1 {\cal V}}(x) \;  {\cal R}_{{\cal V}{\cal V}'}(x,y)
 \;\; \in \;\; {\rm End}(V \otimes {\cal V} \otimes {\cal V}')
  \end{array}
 \ee
 \be
 \lb{Zam23a}
  R_{12}^{-1}(x,y) \;  T_{1 {\cal V}}(x) \;  T_{2 {\cal V}}(y) =
  T_{2{\cal V}}(y)\;  T_{1 {\cal V}}(x) \;  R_{12}^{-1}(x,y)
 \;\; \in \;\; {\rm End}(V \otimes V \otimes {\cal V}) \; ,
\ee
 where we
 denote by ${\cal V}'$ the second copy of the vector space ${\cal V}$,
 the numbers $1,2$ numerate vector spaces $V$, and the matrix
 $R_{12}(x,y)\in {\rm End}(V \otimes V)$, as well as the matrix
 ${\cal R}(x,y) \in {\rm End}({\cal V} \otimes {\cal V}')$, satisfy
 the Yang-Baxter equation (\ref{Zam04}). Consider the transfer-matrix
 \be
 \lb{Zam24}
 \tau(z_1, \dots , z_n) = {\rm Tr}_{\cal V}
 \Bigl( T_{n\cal V}(z_n) \cdots T_{2\cal V}(z_2) T_{1\cal V}(z_1) \; \rho_{\cal V} \Bigr) \; ,
 \ee
 where the operator $\rho_{\cal V} \in {\rm End}({\cal V})$ is such that
 \be
 \lb{Zam25}
 {\cal R}_{{\cal V}{\cal V}'}(x,y) \rho_{\cal V} \rho_{\cal V'} =
 \rho_{\cal V} \rho_{\cal V'} {\cal R}_{{\cal V}{\cal V}'}(x,y) \; .
 \ee
 Then we have   \\

\noindent
 {\bf Proposition~ \ref{sec4}.3.}~~ {\it Transfer-matrices $\tau(z_1, \dots , z_n)$  and
 $\tau(z_1', \dots , z_n')$, defined in (\ref{Zam24}), are commutative
 generating functions
 \be
 \lb{Zam25rt}
 [ \tau(z_1, \dots , z_n) , \; \tau(z_1', \dots , z_n')] = 0 \; ,
 \ee
 if parameters $(z_1, \dots , z_n)$, $(z_1', \dots , z_n')$ and
 the matrix ${\cal R}(x,y)$ are such that
 \be
 \lb{Zam25r}
 {\cal R}(z_n,z_n') = {\cal R}(z_k,z_k') \;\;\;\; \forall k=1,2,\dots,n-1 \; .
 \ee}
 {\bf Proof.} Let ${\cal V}'$ be the second copy of the space ${\cal V}$. Then we have
 $$
 \begin{array}{c}
 \tau(z_1, \dots , z_n)  \; \tau(z_1', \dots , z_n') =
 {\rm Tr}_{{\cal V}{\cal V}'} \Bigl( T_{n\cal V}(z_n) T_{n\cal V '}(z_n') \cdots
 T_{1\cal V}(z_1) T_{1\cal V'}(z_1') \; \rho_{\cal V} \rho_{\cal V'} \Bigr) = \\ [0.3cm]
 =
 {\rm Tr}_{{\cal V}{\cal V}'} \Bigl( {\cal R}^{-1}_{\cal V V'}(z_n,z_n') \cdot
 T_{n\cal V '}(z_n') T_{n\cal V}(z_n)  \cdots
 T_{1\cal V'}(z_1') T_{1\cal V}(z_1)  \cdot
  {\cal R}_{\cal V V'}(z_1,z_1') \; \rho_{\cal V'} \rho_{\cal V}  \Bigr) = \\ [0.3cm]
  = {\rm Tr}_{{\cal V}{\cal V}'} \Bigl(T_{n\cal V '}(z_n') T_{n\cal V}(z_n)  \cdots
 T_{1\cal V'}(z_1')  T_{1\cal V}(z_1) \; \rho_{\cal V'} \rho_{\cal V}  \Bigr) =
 \tau(z_1', \dots , z_n') \;  \tau(z_1, \dots , z_n) \; ,
 \end{array}
 $$
 where ${\rm Tr}_{{\cal V}{\cal V}'}={\rm Tr}_{{\cal V}}{\rm Tr}_{{\cal V}'}$
  and we have used relations (\ref{Zam23}), (\ref{Zam25}). \hfill \qed

 \vspace{0.3cm}

 Note that for the rational (or trigonometric) $R$-matrices, when we have
 $R(x,y) = R(x-y)$ (or $R(x,y) = R(x/y)$), relation
 (\ref{Zam25r}) is fulfilled for the choice $z_k - z_k' = x - y$
 (or $z_k / z_k' = x/y$) for all $k$, where $x$ and $y$ are two fixed parameters.
 For example, in the trigonometric case the commutative transfer-matrix can be
 taken in the form $\tau(x; z_1, \dots , z_n) = \tau(x \, z_1, \dots , x \, z_n)$ and commutativity
 condition (\ref{Zam25rt}) is written as
 $$
 [\tau(x; z_1, \dots , z_n), \; \tau(y; z_1, \dots , z_n)] = 0 \; .
 $$

 Now,  in addition to the relation (\ref{Zam25}), we require that the operator $\rho_{\cal V}$
  satisfies (cf. (\ref{Zam20})):
 \be
 \lb{Zam26}
 \begin{array}{c}
 T_{1\cal V}(x)  \; \rho_{\cal V} \; Q_1  =
 \rho_{\cal V} \; Q_1 \; T_{1\cal V}(\tilde{\bar{x}})  \; , \;\;\;\;
 Q \in {\rm End}(V)
  \end{array}
 \ee
 where for the invertible matrix $Q$ we have (cf. \ref{Zam21})
 $$
 R_{12}(x,y) Q_1 Q_2 = Q_1 Q_2 R_{12}(\tilde{\bar{x}},\tilde{\bar{y}})  \; .
 $$
 Equation (\ref{Zam26}) serves twisted periodic boundary conditions
 of the type (\ref{percond}) for the
 transfer-matrix  (\ref{Zam24}).

 At the end of this Section we formulate the following statement.

  \vspace{0.3cm}

\noindent
  {\bf Proposition~ \ref{sec4}.4.}~~ {\it Flat connections (\ref{Zam22f}) commute with the
   transfer-matrix (\ref{Zam24})
 \be
 \lb{Zam33}
 {\sf A}_k(z_1,\dots,z_n) \; \tau(z_1, \dots , z_n) =
 \tau(z_1, \dots , z_n) \; {\sf A}_k(z_1,\dots,z_n)  \; .
 \ee}
{\bf Proof.} Take the operator
 $T_{k\cal V}(z_k)$ (in the right hand side of (\ref{Zam24})) and use the same procedure as in Remark 2.
 for the cyclic moving of $T_{k\cal V}(z_k)$.
After direct calculations with the use of the relations (\ref{Zam26}),
(\ref{Zam23a}) and identity
$$
\tau(z_1, \dots , \tilde{\bar{z}}_k, \dots , z_n) =
D_{z_k} \; \tau(z_1, \dots , z_k, \dots , z_n) \; D_{z_k}^{-1} \; ,
$$
we deduce relation (\ref{Zam33}).
 \hfill \qed

 \vspace{0.2cm}

 \noindent
 {\bf Consequence.} By using the statement of the Proposition \ref{sec4}.4 we deduce
 the following result. Let $\Psi(z_n,\dots,z_1)$ be any solution of the
 periodic quantum Knizhnik-Zamolodchikov equation (\ref{Zam22}). Then, the vector
 $$
 \Psi'(z_n,\dots,z_1) = \tau(z_1, \dots , z_n) \cdot \Psi(z_n,\dots,z_1) \; ,
 $$
 is also a solution of the
 periodic quantum Knizhnik-Zamolodchikov equation (\ref{Zam22}).

 \section{Flat connections
 for Birman--Murakami--Wenzl algebra.\label{aBMW}}
\setcounter{equation}0



\paragraph{1. Birman--Murakami--Wenzl algebra. Definition and basic relations.} \hspace{-.2cm}
${}$

The {\it Birman--Murakami--Wenzl algebra}
$BMW_n(q,\nu)$ was defined in \cite{BirWen}, \cite{Mur00}
and \cite{Mur0}. It is generated over $\mathbb{C}$ by invertible
elements ${\sf T}_1,\dots,{\sf T}_{n-1}$ with the following defining relations
\be\lb{bmw01}
 {\sf T}_i\, {\sf T}_{i+1}\, {\sf T}_i={\sf T}_{i+1}\, {\sf T}_i\, {\sf T}_{i+1}\; ,\;\;\;
  {\sf T}_i\, {\sf T}_{j}={\sf T}_j\, {\sf T}_{i}\;\;\; {\rm for}\;\; |i-j|>1\; ,
  \ee
\be
\lb{bmw1}
 \kappa_i\, {\sf T}_i={\sf T}_i\,\kappa_i=\nu\,\kappa_i\; ,\phantom{{\sf T}_j}
  \ee
 \be
 \lb{bmw2}
 \kappa_i\, {\sf T}_{i-1}^{\varepsilon}\,\kappa_i=\nu^{-\varepsilon}\,\kappa_i\;\; ,\;\;
\kappa_i\, {\sf T}_{i+1}^{\varepsilon}\,\kappa_i=\nu^{-\varepsilon}\,\kappa_i\ \ \text{with}\ \varepsilon=\pm 1\; ,
 \ee
where
 \be
 \lb{bmw3}\kappa_i:=1-\frac{{\sf T}_i-{\sf T}_i^{-1}}{q -q^{-1}}\; .
 \ee
Here $q$ and $\nu$ are complex parameters of the algebra which we assume generic in the sequel; in particular, the
definition
(\ref{bmw3}) makes sense, the denominator in the right hand side does not vanish.
  Note that the algebra $BMW_n(q,\nu)$ with defining relations
  (\ref{bmw01})--(\ref{bmw2}) possesses the automorphism $\rho_2({\sf T}_i) = {\sf T}_{n-i}$
  (cf. (\ref{autom})).

\vskip .2cm
The quotient of the algebra $BMW_n(q,\nu)$ by the ideal generated by the elements $\kappa_1,\dots,\kappa_n$
(in fact, this ideal is generated by any one of these elements, say, $\kappa_1$) is isomorphic to the Hecke algebra $H_n(q)$. It is also well-known that the
braid group ${\cal{B}}_n$ (of type $A$) embeds in the $BMW_n$ algebra ~~${\cal{B}}_n \hookrightarrow BMW_n$.~~
We shall often omit the parameters in the notation for the algebras and write simply $BMW_n$
and $H_n$.

\vskip .2cm
Let
\be\label{lammu}\mu =\displaystyle{\frac{q -q^{-1}+ \nu^{-1} - \nu}{q -q^{-1}}}=
\displaystyle{\frac{(q^{-1}+\nu)(q-\nu)}{\nu(q -q^{-1})}}\ .\ee
The following relations can be derived from  (\ref{bmw01})--(\ref{bmw2}):
\begin{equation}\lb{bmw4}\kappa_i^2 =  \mu \, \kappa_i
\end{equation}
then, with $\varepsilon =\pm 1$,
\begin{eqnarray}
\lb{bmw5}\kappa_i\, {\sf T}_{i+\varepsilon}\, {\sf T}_i\!\!  &=&\!\!  {\sf T}_{i+\varepsilon}\, {\sf T}_i \,\kappa_{i+\varepsilon} \;
,\\[.5em]
\lb{bmw8}\kappa_i\,\kappa_{i+\varepsilon}\,\kappa_i\!\!  &=&\! \! \kappa_{i}\; ,\\[.5em]
\lb{bmw9}\bigl({\sf T}_i-\!(q -q^{-1})\bigr)\kappa_{i+\varepsilon}\bigl({\sf T}_i-\!(q -q^{-1})\bigr)\!\!   &=&\!\!
\bigl( {\sf T}_{i+\varepsilon}-\!(q -q^{-1})\bigr)\kappa_{i}\bigl({\sf T}_{i+\varepsilon} -\!(q -q^{-1})\bigr)\, ,\\[.5em]
\lb{bmw8a}{\sf T}_{i+\varepsilon}\,\kappa_i\, {\sf T}_{i +\varepsilon}\!\!  &=&\!\!  {\sf T}^{-1}_i\,\kappa_{i+\varepsilon}\, {\sf T}^{-1}_i\;
,
\end{eqnarray}
and
\begin{eqnarray}
\lb{bmw6}\kappa_i\, {\sf T}_{i+\varepsilon}\, {\sf T}_{i}&=&\kappa_i\,\kappa_{i+\varepsilon}\; ,\\[.5em]
\lb{bmw7}\kappa_i\, {\sf T}^{-1}_{i+\varepsilon}\, {\sf T}^{-1}_{i}&=&\kappa_i\,\kappa_{i+\varepsilon}\; ,\\[.5em]
\lb{bmw8b}\kappa_{i+\varepsilon}\,\kappa_i\,\bigl({\sf T}_{i+\varepsilon}-(q -q^{-1})\bigr) &=&
\kappa_{i+\varepsilon}\,\bigl( {\sf T}_i-(q -q^{-1})\bigr)\; ,\end{eqnarray}
together with their images under the anti-automorphism $\rho_a$ of the algebra $BMW_n$ defined on the generators by
\be\lb{aau}
 \rho_a({\sf T}_i)={\sf T}_i\; ,\;\;\;
 \rho_a({\sf T}_i \, {\sf T}_k)={\sf T}_k \, {\sf T}_i\; ,\;\;\;
  \rho_a({\sf T}_i \, {\sf T}_j \, {\sf T}_k)=
  {\sf T}_k \, {\sf T}_j\, {\sf T}_i \; ,\;\;\; \dots .
 \ee

\paragraph{2. Baxterized elements.}

${}$

The {\it baxterized elements} $T_i(u,v)$ are defined by
 \be
\lb{a00}
 T_i(u,v):={\sf T}_i+\frac{q-q^{-1}}{v/u-1}+
\frac{q-q^{-1}}{1+\nu^{-1}qv/ u} \kappa_i \equiv T_i(u/v) \; ,
 \ee
see  \cite{CGX}, \cite{Isa}, \cite{Jones} and \cite{Mur}. They are rational functions in complex variables $u$ and $v$
which
are called {\it spectral variables}. The elements $T_i(u,v)$ depend on the ratio of the spectral parameters; for us it is more convenient to
 have both spectral variables in the notation (\ref{a00}) for the baxterized element.
However for brevity we shall
denote sometimes the baxterized elements by $T_i(u/v) \equiv T_i(u,v)$ (with one argument only).

\vskip .2cm
The baxterized elements satisfy the braid relation of the form
 \be\lb{ybe0}
 T_i(u_{2},u_{3})T_{i+1}(u_1,u_{3})T_i(u_1,u_{2})=
  T_{i+1}(u_1,u_{2})T_i(u_1,u_{3})T_{i+1}(u_{2},u_{3})\; .
  \ee
The inverses of the baxterized elements are given by
\be\lb{aa}T_i(v,u)^{-1}=T_i(u,v)\, f(u,v)\; ,\ee
where
\be f(u,v)=\frac{(u-v)^2}{(u-q^2 v)(u-q^{-2}v)}=f(v,u)\ .\ee

\paragraph{3. Jucys--Murphy elements.}

${}$

The Jucys--Murphy elements of the algebra $B\!M\!W_n$ are defined by
\be\lb{jme}y_1=1 \;, \;\;\; y_{k+1}={\sf T}_k\dots {\sf T}_2\, {\sf T}_1^2\, {\sf T}_2\dots {\sf T}_k\ , \ \ k=1, \dots ,n-1\; .\ee
The elements $y_1,\dots,y_n$ pairwise commute and satisfy the identities
\be
\lb{yyk}
\kappa_j\, y_{j+1}\, y_j=y_j\, y_{j+1}\,\kappa_j=\nu^2\,\kappa_j\; .
\ee

The Jucys--Murphy elements were originally used for constructing idempotents for the symmetric groups in \cite{Ju}, \cite{Mu}.
Analogues of the Jucys--Murphy elements can be defined for a number of important algebras related to the symmetric
group rings (e.g., the Hecke and Brauer algebras);
they turn out to generate maximal commutative subalgebras in these rings (see \cite{IMO}, \cite{IsOg3}, \cite{OPdA}, \cite{OV}
and references therein).
The commutative subalgebra, generated by the Jucys--Murphy elements $y_1,\dots,y_n$, of the generic
algebra $B\!M\!W_n$  is maximal as well; it follows from the results in \cite{IsOg},\cite{LeRa}.


 \paragraph{4. Affine BMW algebras  of type $C$} (see, e.g., \cite{IsOg} and
references therein).

${}$

Affine Birman-Murakami-Wenzl algebras
 $BMW_{n}(C)$ of type $C$ are extensions of the algebras
$BMW_{n}$. The algebra  $BMW_{n}(C)$ is
generated by the elements $\{{\sf T}_1,\dots,{\sf T}_{n-1}\}$ with relations (\ref{bmw01}),
(\ref{bmw1}), (\ref{bmw2}), (\ref{bmw3}) and the affine element ${\sf T}_0 = y_1 \neq 1$ which satisfies
\be
\lb{bmw02}
\begin{array}{c}
{\sf T}_1 \, {\sf T}_0 \, {\sf T}_1 \, {\sf T}_0 = {\sf T}_0 \, {\sf T}_1 \, {\sf T}_0 \, {\sf T}_1 \; , \;\;\;
[{\sf T}_k, \, {\sf T}_0]=0 \;\;\; \mathrm{for} \;\; k > 1 \; ,
\\ [0.2cm]
\kappa_1 \; {\sf T}_0 \, {\sf T}_1 \, {\sf T}_0 \, {\sf T}_1 = \hat{z} \, \kappa_1 =
{\sf T}_1 \, {\sf T}_0 \, {\sf T}_1 \,  {\sf T}_0 \; \kappa_1  \; ,
\end{array}
\ee
\be
 \lb{bmw02z}
\kappa_1 \; {\sf T}_0^k \; \kappa_1 = \hat{z}^{(k)} \kappa_1\ \ ,
\;\;\; k=1,2,3,\dots  \; ,
 \ee
where $\hat{z}$, $\hat{z}^{(k)}$ are central elements. Initially, the affine
version of the Brauer algebras (which are the special limit $q \to 1$
of $BMW_{n}(C)$),  was introduced by M. Nazarov \cite{Nazar}.
Note that the central elements
$\{ \hat{z}^{(k)} \}$ generate an infinite dimensional
abelian subalgebra in $BMW_{n+1}(C)$ and we denote this subalgebra as $BMW_{0}(C)$.

Consider the set of affine elements
(cf. with elements $a_i$
in (\ref{jucys1}))
$$
 y_1 = {\sf T}_0 \; , \;\;\;
 y_{k+1} = {\sf T}_k \, y_k \, {\sf T}_k \; , \;\;\; k=1,2,\dots ,n-1 \; .
$$
The elements $y_k$ $(k=1,2,\dots,n)$ generate a commutative subalgebra $Y_{n}$ in
$BMW_{n}(C)$.

\noindent
{\bf Proposition \ref{aBMW}.1} \cite{IsOg2},\cite{IMO2} {\it For the affine
 BMW algebra the element
\be\lb{rea15}
 L_j(u)=\frac{u-y_{j}}{cu\, y_{j} - 1}
 \; , \;\;\;\;\;
c=- \nu \, q^{-1} \, \hat{z}^{-1} \; ,
\ee
is the baxterized  solution of the reflection equation
 \be
 \lb{rea14}
T_{j}(u,v) \, L_{j}(u) \, T_{j}(v,\bar{u} ) \, L_{j}(v)
= L_{j}(v) \, T_{j}(u,\bar{v} ) \, L_{j}(u) \,
T_{j}(u,v)  \; ,
\;\;\;\;\; (j=1,\dots,n-1) \; ,
 \ee
where $\bar{u} = 1/(c \, u)$.} \\
{\bf Proof.} The formula (\ref{rea14}) is checked by direct calculations. \hfill \qed

 Since we have $T_j\bigl(u,v \bigr) = T_j\bigl(\bar{v},\bar{u} \bigr)$, the equation (\ref{rea14}) is a
 realization of the reflection equation (\ref{Zam10})
 if we identify
 $$
 L_j(v) \to K_j(v) \; , \;\;\;
T_j\bigl(u,v \bigr) \to  \hat{R}_j(u,v)
\; , \;\;\;  \bar{x} = \sigma(x) = \frac{1}{c \, x} \; .
 $$

 \paragraph{5. Affine BMW algebras  of type $C^{(1)}$.}

${}$

The algebra  $BMW_{n}(C^{(1)})$ is generated by
 the elements $\{T_0,T_1, \dots,T_n \}$ and is associated to the
 Coxeter graph (\ref{C1}) of type $C^{(1)}$. The algebra
 $BMW_{n}(C^{(1)})$ is an extension of the
 affine algebra $BMW_{n}(C)$ (we add new generator $T_n$). We require that
 the algebra $BMW_{n}(C^{(1)})$ possesses
 the automorphism $\rho_2$ which is defined in (\ref{autom}). Thus, applying
 automorphism $\rho_2$
 to the relations (\ref{bmw02}), we obtain relations
 for the affine element $T_n$ in the form
 \be
\lb{bmw02n}
\begin{array}{c}
{\sf T}_{n-1} \, {\sf T}_n \, {\sf T}_{n-1} \, {\sf T}_n =
{\sf T}_n \, {\sf T}_{n-1} \, {\sf T}_n \, {\sf T}_{n-1} \; , \;\;\;
[{\sf T}_k, \, {\sf T}_n]=0 \;\;\; \mathrm{for} \;\; k < n-1 \; ,
\\ [0.2cm]
\kappa_{n-1} \; {\sf T}_n \, {\sf T}_{n-1} \, {\sf T}_n \, {\sf T}_{n-1}
= \hat{z}' \, \kappa_{n-1} =
{\sf T}_{n-1} \, {\sf T}_n \, {\sf T}_{n-1} \,  {\sf T}_n \; \kappa_{n-1}  \; ,
\\ [0.2cm]
\kappa_{n-1} \; {\sf T}_n^k \; \kappa_{n-1} =
 \hat{z}^{\prime (k)} \kappa_{n-1}\ \ ,
\;\;\; k=1,2,3,\dots  \; .
\end{array}
\ee
where $\hat{z}' = \rho_2(\hat{z})$, $\hat{z}^{\prime(k)} = \rho_2(\hat{z}^{(k)})$
(as well as  $\hat{z}$, $\hat{z}^{(k)}$)
are the central elements in the algebra $BMW_{n}(C^{(1)})$.

Consider the set of affine elements (cf. with elements $b_i$
in (\ref{jucys1}))
$$
 \bar{y}_n = {\sf T}_n = \rho_2({\sf T}_0) \; , \;\;\;
 \bar{y}_{k-1} = {\sf T}_{k-1} \, \bar{y}_k \, {\sf T}_{k-1} = \rho_2(y_{n-k+2}) \; , \;\;\;
 k=2,\dots ,n \; .
$$
The elements $\bar{y}_k$ $(k=1,\dots,n)$ generate a commutative subalgebra $\bar{Y}_{n}$ in
$BMW_{n}(C^{(1)})$.

Since the element (\ref{rea15})
 is a solution of the reflection elution (\ref{rea14}),
 the element
\be\lb{rea151}
 \bar{L}_j(u)=\frac{u-\bar{y}_{j}}{c' \, u\, \bar{y}_{j} - 1} = \rho_2(L_{n-j+1})
 \; \in \; BMW_{n}(C^{(1)})
 \; , \;\;\;\;\;
c' =-\nu \, q^{-1} \, \hat{z}^{\prime -1} = \rho_2(c) \; ,
\ee
is the baxterized  solution of the
 dual reflection equation
 which is obtained as the image of (\ref{rea14}) under the automorphism $\rho_2$
\be\lb{rea14d}
T_{j}(u,v) \, \bar{L}_{j+1}(u) \, T_{j}(v,\tilde{u} ) \, \bar{L}_{j+1}(v)
= \bar{L}_{j+1}(v) \, T_{j}(u,\tilde{v} ) \, \bar{L}_{j+1}(u) \,
T_{j}(u,v)  \; ,
\;\;\;\;\; (j=1,\dots,n-1) \; ,
\ee
where $\tilde{u} = 1/(u \, c^{\prime})$. Taking into account relations (\ref{aa})
and identities
 $$
 T_j(\tilde{v},u) = T_j(\tilde{u},v) \; , \;\;\;
 \bar{L}_{j}(\tilde{u}) = \frac{1}{c'} \;  \bar{L}_{j}(u)^{-1}
 $$
 we write (\ref{rea14d}) in the form
\be\lb{rea14dd}
T_{j}(v,u) \, \bar{L}_{j+1}(\tilde{u}) \, T_{j}(\tilde{v},u ) \, \bar{L}_{j+1}(\tilde{v})
= \bar{L}_{j+1}(\tilde{v}) \, T_{j}(\tilde{u},v ) \, \bar{L}_{j+1}(\tilde{u}) \,
T_{j}(\tilde{u},\tilde{v})  \; ,
\;\;\;\;\; (j=1,\dots,n-1) \; .
\ee
The equation (\ref{rea14dd}) can be represented
 as the reflection equation (\ref{Zam10a})
 if we identify
 $$
 \bar{L}_j(\tilde{v}) \to K_j(v) \; , \;\;\;
T_j\bigl(u,v \bigr) \to  \hat{R}_j(u,v)
\; , \;\;\;  \tilde{x} = \bar{\sigma}(x) = \frac{1}{c' \, x} \; .
 $$

 \paragraph{6. Embedding of the braid group $B_n(C^{(1)})$
 into the algebra $BMW_{n}(C^{(1)})$.}

${}$

 Let $\{z_1,\dots,z_n \}$ be a set of spectral parameters. Consider the
 Weyl group generated by the operators $s_i$
 (see (\ref{Affbg1})) and the elements
 $$
 T_i(z_i,z_{i+1}) \; , \;\; L_1(z_1) \; , \;\; \bar{L}_n(z_n) \;\; \in \;\; BMW_{n}(C^{(1)}) \; .
 $$
Then we have the following statement.

 \noindent
 {\bf Proposition \ref{aBMW}.2}
 The map $\rho_{b}$ of the
 affine braid group $B_n(C^{(1)})$ into  $BMW_{n}(C^{(1)})$
 defined as (cf. (\ref{Zam10b}))
 \be
 \lb{Zam10r}
 \begin{array}{c}
 \rho_{b}(T_i) \; = \;  s_i \, T_{i}(z_i,z_{i+1})  \;\;\;\; (i=1,\dots,n-1) \; , \;\;\;
 \rho_{b}(T_0) \; = \;   L_1(z_1) \, s_0  \; , \;\;\;
 \rho_{b}(T_n) \; = \;  s_n \, \overline{L}_n(z_n)   \; ,
 \end{array}
 \ee
 is the representation of $B_n(C^{(1)})$. \\
 {\bf Proof.}
 One can  directly check that $\rho_{b}(T_i)$ $(i=0,\dots,n)$ satisfy
 defining relations in (\ref{Affbg}), (\ref{Affbg2}) if $T_{k}(z_k,z_{k+1})$,
$L_1(z_1)$ and $\overline{L}_n(z_n)$
  satisfy relations (\ref{ybe0}), (\ref{rea14}) and (\ref{rea14d}), respectively.
  \hfill \qed

  \vspace{0.2cm}

  \noindent
  {\bf Corollary.}  The map $\rho_{c}$ of the
 affine braid group $B_n(C)$ into  $BMW_{n}(C)$
 \be
 \lb{Zam10D}
 \begin{array}{c}
 \rho_{c}(T_i) \; = \;  s_i \, T_{i}(z_i,z_{i+1})  \;\;\;\; (i=1,\dots,n-1) \; , \;\;\;
 \rho_{c}(T_0) \; = \;   L_1(z_1) \, s_0  \; ,
 \end{array}
 \ee
 is the representation of $B_n(C)$.

  \paragraph{7. Flat connections for the algebra $BMW_{n}(C^{(1)})$.}

${}$

  Flat connections for the algebra $BMW_{n}(C^{(1)})$ are defined as
  images $\rho_b(J_i)$ and $\rho_b(\overline{J}_i)$
  of the elements $J_i$ and $\overline{J}_i$ (see (\ref{jucys1}))
  which form the commuting sets
  of elements in affine braid group $B_n(C^{(1)})$.
  The explicit formulas are (cf. (\ref{Zam05f}))
  \be
 \lb{Zam77}
  {\sf A}_k(z_1,\dots,z_n) = \rho_b(J_k) = K_k(z_k ; z_1, \dots , z_{k-1}) \cdot
  \overline{\sf K}_k(\bar{z}_k ; z_{k+1}, \dots , z_{n})  \; ,
  \ee
  where (cf. (\ref{Zam05d}), (\ref{Zam05K}))
   \be
 \lb{Zam05db}
  \begin{array}{c}
 K_{k}(z_k; \, \vec{z}_{(1,k-1)})  = T_{k-1}^{-1}(z_{k-1},z_k) \cdots  T_1^{-1}(z_1,z_k) \; L_1(z_k) \; T_{1}(z_1,\bar{z}_k) \cdots T_{k-1}(z_{k-1},\bar{z}_k) = \\ [0.2cm]
=  \rho_b(T_{k-1}^{-1} \cdots T_1^{-1} \, T_0 \, T_1 \cdots T_{k-1}) \cdot
 (s_{k-1} \cdots s_1 \, s_0 \, s_1 \cdots s_{k-1}) = \rho_b(\bar{a}_k) \cdot {\bf s}(a_k) \; ,
  \end{array}
 \ee
  \be
 \lb{Zam05Kb}
 \begin{array}{c}
 \overline{\sf K}_{k}(\bar{z}_k; \, \vec{z}_{(k+1,n)}) =
\\ [0.3cm]
 = T_{k}(z_{k+1},\bar{z}_k) \cdots T_{n-1}(z_{n},\bar{z}_k) \; \overline{L}_n(\bar{z}_k) \;
 \cdot D_{z_k} \; \cdot
T_{n-1}(z_k ,z_n)  \cdots T_{k}(z_k ,z_{k+1}) = {\sf s}(a_k) \cdot \rho_b(b_k) \; ,
\end{array}
 \ee
 where the finite difference operator $D_{z_k}$ was defined in (\ref{Zam57}) with
 $$
 \tilde{\bar{z}} = \bar{\sigma}(\sigma(z)) = \bar{\sigma}\Bigl(\frac{1}{c\, z}\Bigr)
 = \frac{c}{c'} \; z \; .
 $$
 We stress that $L_j(u) = D_u$ and $\bar{L}_{j+1}(u) = D_u^{-1}$
 (as well as $L_j(u) = 1$ and $\bar{L}_{j+1}(u) = 1$) are
 solutions of the reflection equations (\ref{rea14}) and (\ref{rea14d}), respectively. For example, we can substitute solution $\bar{L}_{j+1}(u) = 1$
 into (\ref{Zam05Kb}) and reduce the flat connection (\ref{Zam77})
 into the form
 \be
 \lb{Zam78}
  {\sf A}_k(z_1,\dots,z_n) =  K_k(z_k ; z_1, \dots , z_{k-1}) \cdot
  \overline{\sf K}_k'(\bar{z}_k ; z_{k+1}, \dots , z_{n})  \; ,
  \ee
  where
  $$
   \begin{array}{c}
 \overline{\sf K}'_{k}(\bar{z}_k; \, \vec{z}_{(k+1,n)}) =
 T_{k}(z_{k+1},\bar{z}_k) \cdots T_{n-1}(z_{n},\bar{z}_k)
 \cdot D_{z_k} \; \cdot
T_{n-1}(z_k ,z_n)  \cdots T_{k}(z_k ,z_{k+1}) \; ,
\end{array}
  $$
\paragraph{8. Braid--Hecke algebra ${\cal{BH}}_n(q,\nu)$.}

${}$

The {\it Braid--Hecke algebra} ${\cal{BH}}_n(q,\nu)$, as far as we know, was
introduced in \cite{Co1},\cite{Co2}. It is generated over $\mathbb{C}$ by the invertible  {\it braid type generators},
$${\sf T_1}, \cdots, {\sf T_{n-1}},$$
 subject the following defining relations

\be\lb{bmw012}
 {\sf T}_i\, {\sf T}_{i+1}\, {\sf T}_i={\sf T}_{i+1}\, {\sf T}_i\, {\sf T}_{i+1}\; ,\;\;\;
  {\sf T}_i\, {\sf T}_{j}={\sf T}_j\, {\sf T}_{i}\;\;\; {\rm for}\;\; |i-j|>1\; ,
  \ee
\be
\lb{bmw12}
 \kappa_i\, {\sf T}_i={\sf T}_i\,\kappa_i=\nu\,\kappa_i\; ,\phantom{{\sf T}_j}
  \ee
 \be
 \lb{bmw22}
 {\sf T}_{i \pm 1}\,{\sf T}_{i} \,\,\kappa_{i \pm 1} - \kappa_{i}\,\kappa_{i \pm 1}
=
{\sf T}_{i}\, {\sf T}_{i \pm 1} \, \kappa_{i} - \kappa_{ i \pm 1}\,\kappa_{i},
\ee
\be
\lb{bmv33}
\kappa_i\,\kappa_{i \pm 1}\,\kappa_{i}- \kappa_{i} = \kappa_{i \pm 1} \,
\kappa_{i}\, \kappa_{i \pm 1} - \kappa_{i \pm 1},
\ee

where
 \be
 \lb{bmw3a}\kappa_i:=1-\frac{{\sf T}_i-{\sf T}_i^{-1}}{q -q^{-1}}\; .
 \ee
Here $q$ and $\nu$ are complex parameters of the algebra which we assume
generic in the sequel; in particular, the  definition  (\ref{bmw3a}) makes
sense, the denominator in the right hand side does not vanish. Note that the
algebra $BMW_n(q,\nu)$ is the quotient of the algebra ${\cal{BH}}_n(1,\nu)$ by the two-sided ideal generated by the {\it tangle relations}
\be
\lb{bmv43}
\kappa_i\,\kappa_{i \pm 1}\,\kappa_{i}- \kappa_{i} = 0, ~~~\kappa_{i \pm 1} \,
\kappa_{i}\, \kappa_{i \pm 1} - \kappa_{i \pm 1} = 0.
\ee

  ~It is easy to see that
$$ ({\sf T}_i -q)({\sf T}_{i}+q^{-1})({\sf T}_{i} - \nu)=0,~~\kappa_{i}^2=
{\frac{(q- \nu)(q^{-1} + \nu)}{\nu (q-q^{-1})}}~~\kappa_{i}.$$
It follows from (\ref{bmv33}), that the elements  $\{ \kappa_1,\ldots, \kappa_{n-1} \}$ generate the Hecke algebra ${\cal{H}}_n(p)$ corresponding to a parameter $p$ such that
$$p+p^{-1} = {\frac{ (q- \nu)(q^{-1} + \nu)}{\nu (q-q^{-1})}}.$$
  Note that the algebra ${\cal{BH}}_n(q,\nu)$ with defining relations
  (\ref{bmw012})--(\ref{bmv33}) possesses the automorphism $\rho_2({\sf T}_i) = {\sf T}_{n-i}$
  (cf. (\ref{autom})). It is well-known that $\dim (BMW_n) = (2n-1) !!$.~~ As
for the algebra ${\cal{BH}}_n(q,\nu)$, it is known \cite{Co2} that  it has finite dimension, but as far as we know, the exact value of its dimension is still unknown. \\
$\bullet$~~~The baxterized elements~ $ \{T_i(u,v),~i=1,\ldots,n-1 \}$~ are defined by (cf (\ref{ybe0}),~~($z:=u/v))$,
$$(\nu+q~z)~ T_i(z) = q~z~{\sf T}_i+\nu~{\sf T}_i^{-1}+(q-q^{-1})~\frac{z(q+\nu)}{z-1}~
 \equiv (\nu+q~z)~T_i(u,v).$$
$\bullet$~~~The Jucys--Murphy elements~ $\{ y_i,~i=1,\ldots,n-1 \}$ ~of the algebra ${\cal{BH}}_n(q,\nu)$ are defined by (\ref{jme}). ~~The $JM$ elements $y_2, ldots, y_{n-1}$ pairwise commute and satisfy the identities
$$\kappa_j~y_{j+1}~y_{j} = y_{j}~y_{j+1}~\kappa_{j},~~~1 \le j < n-1.$$
$\bullet$~~~Affine braid-Hecke algebra  ${\cal{BH}}_n(C)$ of type $C$ is an
 extension of the algebra ${\cal{BH}}_n(q,\nu)$ by the the affine element ${\sf T_{0}}=y_{0} \not= 1$, subject to the set of ``crossing relations'' (\ref{bmw02})~and~(\ref{bmw02z}). One can check that the set of elements
$$y_1= {\sf T}_{0},~~ y_{k+1}:= {\sf T}_{k}~y_k~{\sf T}_{k},~~~k=1,\ldots,n-1 $$
generate a commutative subalgebra in ${\cal{BH}}_n(C)$.

$\bullet$~~(Markov trace, cf \cite{Co2})~~~The family of algebras $\{ {\cal{BH}}_n(q,\nu) \hookrightarrow {\cal{BH}}_{n+1}(q,\nu) \}_{n \ge 1}$ can be provided with (unique) set of homomorphisms
$$  Tr_{n+1} : {\cal{BH}}_{n+1}(q,\nu) \longrightarrow {\cal{BH}}_n(q,\nu) $$
which satisfy the conditions stated in Proposition~\ref{sec6}.2.

\underline{Summarizing}, the all properties of the algebra ${\cal{BH}}_n(q,\nu)$ stated in the item ${\bf 8}$ allow the use of the methods developed in
 Sections~\ref{aBMW} and \ref{sec6},  to construct families of {\it commutative subalgebras}
 in the algebra ${\cal{BH}}_n(q,\nu)$, as well as  ${\cal{BH}}_n(q,\nu)$-valued flat connections and associated $qKZ$ equations. Details will appear.

\vskip .2cm

 \section{Sklyanin's transfer-matrices for affine BMW algebra.\label{sec6}}
\setcounter{equation}0

\vspace{0.2cm}

In this Section, to simplify formulas we make the redefinition
of all spectral parameters $z \to c^{-1/2} z$. In this case the
baxterized element (\ref{a00}) does not changed (since it depends on
the ratio of spectral parameters) and
statement of the Proposition {\bf \ref{aBMW}.1} reads as the following.
For the affine BMW algebra $BMW_{n}(C^{(1)})$ the element
\be\lb{rea15z}
 L_j(c^{-1/2} u)=\frac{c^{-1/2} u-y_{j}}{c^{1/2} u\, y_{j} - 1}
 \equiv y_j(u)
 \;\; , \;\;\;\;\;\; y_j(u) \cdot y_j(u^{-1}) = c^{-1} \; ,
\ee
where $c=- \nu \, q^{-1} \, \hat{z}^{-1}$, is the baxterized  solution of the reflection equation
 \be
 \lb{rea14z}
T_{j}(u/v) \, y_{j}(u) \, T_{j}(v \, u) \, y_{j}(v)
= y_{j}(v) \, T_{j}(v \, u ) \, y_{j}(u) \,
T_{j}(u/v)  \; ,
\;\;\;\;\; (j=1,\dots,n-1) \; .
 \ee

\subsection{Sklyanin's transfer-matrix elements for the algebra $BMW_n(C)$}

In this Subsection we generalize to the BMW algebra case
 results obtained in \cite{Isa07},\cite{IK} for the Hecke algebra case.

{\bf Definition~\ref{sec6}.1}~~ Let  $\vec{z}_{(k)}=(z_1,\dots,z_{k})$ be $k$
parameters and $y_1(x) \in BMW_{n}(C^{(1)})$
is any {\it local}
(i.e., $[y_1(x), T_k]=0,$ $\forall k > 1$) solution
of the reflection equation
(\ref{rea14z}) with $j=1$:
\be
\lb{reflH1}
T_1 \left(x /z \right) \, y_{1}(x) \, T_1(x \, z) \, y_{1}(z) =
y_{1}(z) \, T_1(x \, z) \, y_{1}(x) \,  T_1 \left( x /z \right) \; ,
 \ee
 where solution $y_1(z)$ is
 given in (\ref{rea15z}) for $j=1$.
Define the elements (cf. (\ref{Zam05db}))
\be
\lb{xxz55}
\begin{array}{c}
y_{k}(x;\vec{z}_{(k-1)}) =
T_{k-1}(\frac{x}{z_{k-1}}) \cdots T_2(\frac{x}{z_2}) T_1(\frac{x}{z_1}) \, y_1(x) \,
T_{1}( x z_1) T_2( x z_2) \cdots T_{k-1}( x z_{k-1}) = \\ [0.2cm]
= T_{k-1}(\frac{x}{z_{k-1}}) y_{k-1}(x;\vec{z}_{(k-2)}) T_{k-1}( x z_{k-1}) \; ,
\end{array}
\ee
which we call ``baxterized'' Jucys--Murphy elements.

{\bf Proposition~\ref{sec6}.1}~~ {\it The
``baxterized'' Jucys--Murphy element (\ref{xxz55}) is a solution of the reflection equation
 \be
 \lb{reflH}
T_k \left(x /z \right) \, y_{k}(x;\vec{z}_{(k-1)}) \, T_k( x \, z) \, y_{k}(z;\vec{z}_{(k-1)}) =
y_{k}(z;\vec{z}_{(k-1)}) \, T_k(x \, z) \, y_{k}(x;\vec{z}_{(k-1)}) \,  T_k \left( x /z \right) \; ,
 \ee
}
{\bf Proof.}~~The case $k=1$ of the equation (\ref{reflH})  corresponds to our
assumption that $y_1(x)$ satisfies the equation (\ref{reflH1}). The general
case follows by induction using the definition (\ref{xxz55}) of elements
$y_{k}(x;\vec{z}_{(k-1)}).$ \hfill  \qed

\vspace{0.2cm}

For example, in the case of the affine BMW algebra
$BMW_{n}(C^{(1)})$, one can use the local solution (\ref{rea15z})
 for $j=1$ (recall that $y_1 = {\sf T}_0$):
\be
\lb{soluH}
 y_1(x)= { c^{- 1/2} x - y_1 \over c^{1/2} \, x \, y_1 - 1} =
 { c^{- 1/2} x - {\sf T}_0 \over c^{1/2} \, x \, {\sf T}_0 - 1}
 \;\; , \;\;\;\;\;\;  y_1(1) = - c^{- 1/2} \; .
 \ee

\vspace{0.2cm}

Further we consider only one-boundary affine BMW algebra
 $BMW_{n}(C)$ of type $C$ which is obtained as  the projection ${\sf T}_n=1$  from
 the two-boundary affine BMW algebra $BMW_{n}(C^{(1)})$ (see paragraph 4 in
 Section {\bf \ref{aBMW}}).

Consider the following inclusions of the subalgebras
$BMW_{1}(C) \subset BMW_{2}(C) \subset \dots \subset BMW_{n+1}(C)$:
$$
\{{\sf T}_0; {\sf T}_1, \dots ,{\sf T}_{k-1}\} \in BMW_{k}(C) \subset BMW_{k+1}(C) \ni
\{{\sf T}_0; {\sf T}_1, \dots ,{\sf T}_{k-1},{\sf T}_k \} \; .
$$
For the subalgebras $BMW_{k+1}(C)$
we introduce linear mapping (quantum trace)
$$
{\rm Tr}_{(k+1)}: \;\; BMW_{k+1}(C) \to BMW_{k}(C) \; , \;\;\;\; (k=1,2,\dots,n) \; ,
$$
which is defined by the formula
 \be
 \lb{qtr}
\kappa_{k+1} \, X_{k+1} \, \kappa_{k+1} = \frac{1}{\nu} \, Tr_{({k+1})}(X_{k+1}) \,
\kappa_{k+1} \; , \;\;\; \forall X_{k+1} \in BMW_{k+1}(C) \; .
 \ee
  \vspace{0.2cm}
 \noindent
  {\bf Proposition \ref{sec6}.2} {\it
 For the map $Tr_{(k+1)}$: $BMW_{k+1}(C) \to BMW_{k}(C)$
 we have the following properties $(\forall X_{k},X_k' \; \in \; BMW_{k}(C) \; ,
 \;\; \forall Y_{k+1} \; \in \; BMW_{k+1}(C))$
 \be
 \lb{qtr1}
 \begin{array}{c}
 Tr_{(k+1)}({\sf T}_{k})=1  \; , \;\;
 Tr_{(k+1)}({\sf T}_{k}^{-1})=\nu^2 \; , \;\;
   {\rm Tr}_{(k+1)} ( X_k ) = \nu \, \mu \, X_k  \; , \\[0.1cm]
     Tr_{(k+1)}(\kappa_k)= \nu  \; ,
   \;\; {\rm Tr}_{(1)} ({\sf T}_0^k)= \nu \, \hat{z}^{(k)} \; ,
 \end{array}
 \ee
 \be
 \lb{qtr2}
 Tr_{(k+1)}({\sf T}_{k} \, X_{k} \, {\sf T}_{k}^{-1}) =
 Tr_{(k)}(X_{k})  =
  Tr_{(k+1)}({\sf T}_{k}^{-1} \, X_{k} \, {\sf T}_{k})   \; ,
 \ee
 \be
 \lb{qtr3}
 Tr_{(k+1)}({\sf T}_{k} \, X_{k} \, \kappa_{k}) =
  Tr_{(k+1)}(\kappa_{k} \, X_{k} \, {\sf T}_{k}) \; ,
 \ee
 \be
\lb{map}
\begin{array}{c}
{\rm Tr}_{(k+1)} ( X_k \cdot Y_{k+1} \cdot X'_k )
= X_k \cdot {\rm Tr}_{(k+1)}(Y_{k+1}) \cdot X_k' \;\;  \, , \\[0.1cm]
{\rm Tr}_{(k)}  {\rm Tr}_{(k+1)} ( {\sf T}_k \cdot Y_{k+1} ) = {\rm Tr}_{(k)}  {\rm Tr}_{(k+1)}
( Y_{k+1} \cdot {\sf T}_k)  \; .
\end{array}
\ee}
 {\bf Proof.} Eqs. (\ref{qtr1}) follow from (\ref{bmw2}), (\ref{bmw4}),
 (\ref{bmw8}) and (\ref{bmw02z}). Using (\ref{bmw6}), (\ref{bmw7})
 and (\ref{bmw8}) we have
 $$
 \begin{array}{c}
 \frac{1}{\nu} \, Tr_{(k+1)}({\sf T}_{k} \, X_{k} \, {\sf T}_{k}^{-1})  \kappa_{k+1} =
 \kappa_{k+1} \, {\sf T}_{k} \, {\sf T}_{k+1} \, X_{k} \, {\sf T}_{k+1}^{-1} \, {\sf T}_{k}^{-1} \, \kappa_{k+1}=
  \\ [0.2cm]
 = \kappa_{k+1} \, \kappa_{k} \,  X_{k} \, \kappa_{k} \, \kappa_{k+1}

  = \frac{1}{\nu} \, Tr_{(k)}(X_{k}) \,
   \kappa_{k+1} \, \kappa_{k} \, \kappa_{k+1}  =
    \frac{1}{\nu} \, Tr_{(k)}(X_{k}) \,
   \kappa_{k+1}  \; ,
 \end{array}
 $$
 which is equivalent to the first equality in (\ref{qtr2})
 (second equality in (\ref{qtr2}) can be proved analogously).
 Eq. (\ref{qtr3}) can be proved in the following way
 $$
 \begin{array}{c}
 \kappa_{k+1} \, {\sf T}_{k} \, X_{k} \, \kappa_{k} \, \kappa_{k+1}=
 \kappa_{k+1} \, \kappa_{k} \, {\sf T}_{k+1}^{-1} \, X_{k} \,
  \kappa_{k} \, \kappa_{k+1}=
  \kappa_{k+1} \, \kappa_{k} \,  X_{k} \, {\sf T}_{k+1}^{-1}
  \, \kappa_{k} \, \kappa_{k+1}
  = \\ [0.2cm]
  = \kappa_{k+1} \, \kappa_{k} \,  X_{k} \, {\sf T}_{k} \, \kappa_{k+1}\; .
 \end{array}
 $$
 The first eq. in (\ref{map}) is evident and
the proof of second eq. in (\ref{map}) is the following.
 First of all for any $Y_{k+1}' \in BMW_{k+1}(C)$ we have
  $$
 \kappa_{k+2} \kappa_k  \kappa_{k+1} \, Y_{k+1}' \, \kappa_{k+1} \kappa_k \kappa_{k+2}
 = \frac{1}{\nu} \kappa_{k+2} \kappa_k  \kappa_{k+1}  {\rm Tr}_{(k+1)}(Y_{k+1}')
  \, \kappa_k \kappa_{k+2} =
$$
$$
= \frac{1}{\nu} \kappa_{k+2} \kappa_k    {\rm Tr}_{(k+1)}(Y_{k+1}')
  \, \kappa_k  = \frac{1}{\nu^2} \kappa_{k+2} \kappa_k
   {\rm Tr}_{(k)} {\rm Tr}_{(k+1)}(Y_{k+1}') \; .
$$
 Then, using this equation
 and relations (\ref{bmw6}), (\ref{bmw7}) we obtain
$$
 \begin{array}{c}
 \frac{1}{\nu^2} \kappa_{k+2} \kappa_k
   {\rm Tr}_{(k)} {\rm Tr}_{(k+1)}(Y_{k+1} \, {\sf T}_k) =
 \kappa_{k+2} \kappa_k  \kappa_{k+1} \, Y_{k+1} {\sf T}_k  \, \kappa_{k+1} \kappa_k \kappa_{k+2}
 = \kappa_{k+2} \kappa_k  \kappa_{k+1} \, Y_{k+1}
 {\sf T}_{k+1}^{-1}  \, \kappa_{k+2} \kappa_k    = \\[0.2cm]
 = \kappa_{k+2} \kappa_k  \kappa_{k+1} \, Y_{k+1}
 {\sf T}_{k+2}  \, \kappa_{k+1} \, \kappa_{k+2} \kappa_k  =
 \kappa_{k+2} \kappa_k  \kappa_{k+1} \,  {\sf T}_{k+2} Y_{k+1}
  \, \kappa_{k+1} \, \kappa_{k+2} \kappa_k  = \\[0.2cm]
 =  \kappa_{k+2} \kappa_k  \kappa_{k+1} \,  {\sf T}_{k} Y_{k+1}
  \, \kappa_{k+1} \, \kappa_{k+2} \kappa_k  = \frac{1}{\nu^2} \kappa_{k+2} \kappa_k
   {\rm Tr}_{(k)} {\rm Tr}_{(k+1)}({\sf T}_k \, Y_{k+1}) \; .
 \end{array}
$$
 \qed

Below we use the following identities for baxterized elements (\ref{a00}):
\be
 \lb{bax1}
 \begin{array}{c}
T_n(x)\, T_n(y)= \frac{(q-q^{-1}) (1-xy)}{(1-x)(1-y)} T_n(x \, y) + 1 +
\frac{ (q-q^{-1}) \,
 \nu \, (x y-1) (\nu x y+q^3)}{(\nu y+q)(\nu x+q)(\nu x y+q)}
  \kappa_n \;\; \Rightarrow \\ [0.2cm]
 T_n(x) = T_n(y) + \frac{(q-q^{-1}) (x-y)}{(y-1)(x-1)}
 + \frac{(q-q^{-1})  \nu  q (x-y)}{
 (\nu y+q) (\nu x +q)}\kappa_n \; ,
\end{array}
\ee
 Note that identity (\ref{aa}) is a consequence of the first relation in (\ref{bax1})
 if we substitute there $y=x^{-1}$ and take into account
 $(1-xy) \; T_n(x \, y) \; \stackrel{\; y \to x^{-1}}{\longrightarrow} \;
 (q-q^{-1})$.

 \vspace{0.2cm}

Using the properties (\ref{map}) of the map ${\rm Tr}_{_{(n+1)}}$ and
relations (\ref{bax1}),  one can prove the Lemma.

\vspace{0.1cm}

\noindent
{\bf Lemma~\ref{sec6}.1}~~ {\it For all
$X_k \in BMW_{k}(C)$ and all spectral parameters $x$ and $z$
the following identity is true:
\be
\lb{mapp1}
 \begin{array}{c}
{\rm Tr}_{_{(k+1)}} \Bigl( T_k(x) \cdot X_k \cdot T_k(z) \Bigr)  =
\frac{(q-1/q) (\nu^2 x z -q^2)}{(x \nu +q)(z \nu +q)}
{\rm Tr}_{_{(k+1)}} \Bigl( T_k \cdot X_k \cdot \kappa_k \Bigr) + \\ [0.2cm]
 + \frac{(q^2 \nu^2 x z -\nu^2  x z +x z \nu^2 q^{-2} +q \nu (x+z) +q^2)}{
 (x \nu+q) (z \nu +q)} \,
{\rm Tr}_{_{(k)}} (X_k) -\frac{(q-1/q)(x z \nu^2 -q^2)
 (x z \nu q^2 +x z \nu^2 q +(x +z)\nu -x z \nu +q)}{
((z\nu+q)(x\nu+q)q(z-1)(x-1))} \, X_k \;  ,
\end{array}
\ee
where $T_k(x)$ and $T_k(z)$ are Baxterized elements  (\ref{a00}).
} \\
{\bf Proof.} Direct calculations with the help of properties
(\ref{qtr1}) -- (\ref{map}). \hfill  \qed

\vspace{0.2cm}

{}From eq. (\ref{mapp1}), for $x \, z = q^2 \nu^{-2}$, we obtain the
"crossing-unitarity relation"
 \be \lb{mapp2}
  \Tr _{_{(k+1)}} \Bigl(
T_k(x) \cdot X_k \cdot T_k \left( q^2 \nu^{-2}/ x\right) \Bigr)  =
 \frac{1}{F(x)} \, \Tr _{_{(k)}} (X_k) \;  ,
\ee
where $F(x) =\frac{ (x \nu + q)^2}{(x \nu + q^3)(x\nu + q^{-1})}$.
 Note that identity (\ref{mapp2}) was obtained in \cite{IsOg2} for
 slightly different definition of the baxterized elements (\ref{a00}).

\vspace{0.2cm}

\noindent
{\bf Proposition~\ref{sec6}.3} (see also \cite{IsOg2}, \cite{Isa07}).
 {\it Let $y_{k}(x) \in BMW_{k}(C)$ be any
solution of the RE
(\ref{reflH}). The operators
\be
\lb{tau11}
\tau_{k-1}(x) = \Tr _{_{(k)}} \left( y_{k}(x) \right) \in BMW_{k-1}(C) \; ,
\ee
 form a commutative family of operators
\be
\lb{comfa}
\Bigl[ \tau_{k-1}(x) \, , \; \tau_{k-1}(z) \Bigr]=0 \;\;\; (\forall x,z) \; ,
\ee
in the subalgebra $BMW_{k-1}(C) \subset BMW_{n}(C)$ $(k < n)$.
}

\noindent
{\bf Proof.} Using properties (\ref{qtr2}),
(\ref{map}) and relations (\ref{mapp2}), (\ref{reflH}) we find
$$
\begin{array}{c}
\tau_{k-1}(x) \, \tau_{k-1}(z) = Tr_{_{(k)}} \left( y_{k}(x) \, \tau_{k-1}(z) \right) = \\ [0.2cm]
= F(x \, z) \,
 \Tr _{_{(k)}}  \left( y_{k}(x) \, \Tr _{_{(k+1)}}
  \Bigl( T_k (x \, z) \, y_{k}(z)
 T_k(q^2 (\nu^2 xz)^{-1}) \Bigr)
\right)=  \\ [0.2cm]
= F(x \, z) \,
 \Tr _{_{(k)}} \Tr _{_{(k+1)}} \Bigl( T_k(x/z) \, y_{k}(x) \, T_k (xz) \, y_{k}(z)
 \, T_k^{-1}(x/z) \, T_k (q^2 (\nu^2 xz)^{-1})
\Bigr)= \\ [0.2cm]
= F(x \, z) \,
 \Tr _{_{(k)}} \Tr _{_{(k+1)}} \left(  y_{k}(z) \, T_k (xz) \, y_{k}(x)
 \,  T_k (q^2 (\nu^2 xz)^{-1})  \right)= \\ [0.2cm]
=
 \Tr _{_{(k)}}\left(  y_{k}(z) \, \tau_{k-1}(x)  \right) = \tau_{k-1}(z) \, \tau_{k-1}(x) \; ,
\end{array}
$$
where $F(x)$ was defined in (\ref{mapp2}). \\
\qed

\vspace{0.3cm}

Now we consider the operators (\ref{tau11})
 \be
 \lb{bethe02}
 \tau_n(x;\vec{z}_{(n)}) = \Tr _{_{(n+1)}} \left( y_{n+1}(x;\vec{z}_{(n)})\right) \in BMW_{n}(C) \; ,
 \ee
 where solution $y_{n+1}(x) \in BMW_{n+1}(C)$ of the reflection equation
is taken in the form (\ref{xxz55}).
  We stress that the elements (\ref{bethe02})
   are nothing but the analogs of  Sklyanin's transfer-matrices \cite{Skl}
 and the coefficients in the expansion of $\tau_n(x;\vec{z}_{(n)})$
 over the variable $x$ (for the homogeneous case $z_k=1$) are the Hamiltonians for the
 open Birman-Murakami-Wenzel chain models with nontrivial boundary conditions.

 For example let us redefine all baxterized elements in (\ref{a00})
 \be
 \lb{bax3}
 T_i(x) \; \to \; \tilde{T}_i(x)  = (1-x)  \,  T_i(x) =
 (1-x) \Bigl({\sf T}_i + \frac{(q-q^{-1})x}{x + \nu^{-1} q} \, \kappa_i \Bigr) +
 (q-q^{-1})x \; ,
 \ee
 such that the new elements $\tilde{T}_i(x)$ satisfies
 conditions
 $$
 \left. \tilde{T}_i(x) \right|_{x=1} = (q-q^{-1}) \;\; , \;\;\;\;\;
 \left. \partial_x \, \tilde{T}_i(x) \right|_{x=1} = -{\sf T}_i
  - \frac{(q-q^{-1})}{1+ \nu^{-1}q} \, \kappa_i + (q-q^{-1}) \; .
 $$
 \be
 \lb{bax4}
 \tilde{T}_i(u/v)  \; \tilde{T}_i(v/u)  =
 \frac{(v q^2 -u)  \,  (u q^2  - v)}{q^2 u v} \; .
 \ee
 Now we respectively redefine the Sklyanin's transfer-matrix element (\ref{bethe02})
 as the following
 \be
 \lb{bethe02z}
 \tilde{\tau}_n(x;\vec{z}_{(n)}) = \prod_{i=1}^n \bigl( (1-x/z_i)
 \left(1-x z_i \right) \bigr)
 \Tr _{_{(n+1)}} \left( y_{n+1}(x;\vec{z}_{(n)})\right)  =
 \Tr _{_{(n+1)}} \left( \tilde{y}_{n+1}(x;\vec{z}_{(n)})\right)  \; ,
 \ee
 where $\tilde{y}_{n+1}(x;\vec{z}_{(n)})$ is given by (\ref{xxz55})
 with substitution $T_i(x)  \to  \tilde{T}_i(x)$ and $k \to n+1$. I.e.,
 we have
 \be
\lb{xxz55w}
\begin{array}{c}
\tilde{y}_{k}(x;\vec{z}_{(k-1)}) =
\tilde{T}_{k-1}(\frac{x}{z_{k-1}}) \cdots \tilde{T}_2(\frac{x}{z_2}) \tilde{T}_1(\frac{x}{z_1}) \, y_1(x) \,
\tilde{T}_{1}( x z_1) \tilde{T}_2( x z_2) \cdots \tilde{T}_{k-1}( x z_{k-1}) = \\ [0.2cm]
= \tilde{T}_{k-1}(\frac{x}{z_{k-1}}) \tilde{y}_{k-1}(x;\vec{z}_{(k-2)})
 \tilde{T}_{k-1}( x z_{k-1}) \; .
\end{array}
\ee
Using "unitarity condition" (\ref{bax4}) we represent
baxterized Jucys-Murphy elements (\ref{xxz55w})
in the form
 \be
\lb{xxz55y}
\begin{array}{c}
\tilde{y}_{k}(x;\vec{z}_{(k-1)}) = \Bigl( \prod\limits_{i=1}^{k-1}
\frac{(x q^2 -z_i)  \,  (z_i q^2  - x)}{q^2 x z_i} \Bigr)
\tilde{y}_{k}'(x;\vec{z}_{(k-1)}) \; ,
 \\ [0.3cm]
\tilde{y}_{k}'(x;\vec{z}_{(k-1)})  \equiv
\tilde{T}_{k-1}^{-1}(\frac{z_{k-1}}{x}) \cdots  \tilde{T}_1^{-1}(\frac{z_1}{x}) \, y_1(x) \,
\tilde{T}_{1}( x z_1) \tilde{T}_2( x z_2) \cdots \tilde{T}_{k-1}( x z_{k-1})
 \; .
\end{array}
\ee
We will use this form below.

  Then,
 for the homogeneous case $z_i=1$ $(\forall i)$,
 we consider the coefficient
 $$
 \frac{c^{1/2} (q-q^{-1})^{1-2n}}{2\nu \mu }  \Bigl(
  \left. \partial_x \, \tilde{\tau}_n(x;z_i=1)\; \right|_{x=1} \Bigr) =
  \sum_{i=1}^{n-1} \Bigl( {\sf T}_i
  + \frac{(q-q^{-1})}{1+ \nu^{-1}q} \, \kappa_i \Bigr)
  + \frac{(q-q^{-1})}{2} \frac{c {\sf T}_0^2 -1}{(c^{1/2} {\sf T}_0 -1)^2} + {\rm constant} \; ,
 $$
 in the expansion of
 the generating function $\tilde{\tau}_n(x;z_i=1)$ for commutative elements.
 This coefficient gives (up to an additional constant) the element
 $$
 {\cal H} = \frac{(q-q^{-1})}{2} \frac{c {\sf T}_0^2 -1}{(c^{1/2} {\sf T}_0 -1)^2} +
 \sum_{i=1}^{n-1} \Bigl( {\sf T}_i
  + \frac{(q-q^{-1})}{1+ \nu^{-1}q} \, \kappa_i \Bigr) \;\; \in \;\;
  BMW_n(C) \; ,
 $$
 being  the local
 Hamiltonian
 for the open BMW chain model with nontrivial boundary
 condition for the first site of the chain.

 Consider the expansion of $\tau_n(x;\vec{z}_{(n)})$
 over $x$ for the inhomogeneous case:
 \be
 \lb{bethe01}
 \tau_n(x;\vec{z}_{(n)})  = \sum_{k=-\infty}^\infty \Phi_k(\vec{z}_{(n)})  \, x^k \;  \in BMW_{n}(C)  \; .
 \ee
 According to the
 Proposition \ref{sec6}.3, for fixed parameters
 $\vec{z}_{(n)}=(z_{1},\dots,z_n)$, the elements $\Phi_k(\vec{z}_{(n)})$
   generate a commutative subalgebra
 $\hat{\cal B}_n(\vec{z}_{(n)}) \subset BMW_{n}(C) $. These elements are interpreted as
 Hamiltonians for the inhomogeneous
 open Hecke chain models. Following \cite{MTV}
 we call the subalgebras $\hat{\cal B}_n(\vec{z}_{(n)})$ as Bethe subalgebras of the affine algebra $BMW_{n}(C)$.

 \subsection{Bethe subalgebras for affine BMW algebra and q-KZ connections}

In this Section and below we will use the normalized baxterized elements (\ref{bax3}):
 $\widetilde{T}_k(x) = (1-x) T_k(x)$.
 Consider the transfer-matrix operator (\ref{bethe02z})
 and fix the spectral parameter as $x=z_k$, where $1 \leq k \leq n$
 (analogous results can be obtained if instead we take $x=z_k^{-1}$).
 In view of relation $\left. T_k(x/z_k)\right|_{x=z_k} = (q-q^{-1})$
 we deduce for the transfer-matrix operator (\ref{bethe02z})
 $$
 \begin{array}{c}
 B_k(\vec{z}) =  \frac{1}{(q-q^{-1})} \tau_{n}(z_k ;\vec{z}_{(n)}) =
\Tr _{_{\!\! (n+1)}} \! \Bigl( \!
\underline{\tilde{T}_{n}(\frac{z_k}{z_{n}}) \cdots \tilde{T}_{k+1}(\frac{z_k}{z_{k+1}})}
\tilde{T}_{k-1}(\frac{z_k}{z_{k-1}}) \cdots \tilde{T}_1(\frac{z_k}{z_1}) \, y_1(z_k) \cdot \\ [0.3cm]
\cdot
\tilde{T}_{1}(z_k z_1) \cdots \tilde{T}_{k-1} (z_k z_{k-1}) \tilde{T}_k(z_k^2)
 \tilde{T}_{k+1} (z_k z_{k+1}) \cdots \tilde{T}_{n}(z_k z_{n})\! \Bigr) = \\ [0.3cm]
\Tr _{_{\!\! (n+1)}} \! \Bigl( \!
\tilde{T}_{k-1}(\frac{z_k}{z_{k-1}}) \cdots \tilde{T}_1(\frac{z_k}{z_1}) \, y_1(z_k) \cdot \\ [0.3cm]
\cdot \tilde{T}_{1}(z_k z_1) \cdots \tilde{T}_{k-1} (z_k z_{k-1}) \cdot
\underline{\tilde{T}_{n}(\frac{z_k}{z_{n}}) \cdots \tilde{T}_{k+1}(\frac{z_k}{z_{k+1}})}  \tilde{T}_k(z_k^2)
 \tilde{T}_{k+1} (z_k z_{k+1}) \cdots \tilde{T}_{n}(z_k z_{n})\! \Bigr) =
 \end{array}
 $$
 \be
 \lb{qkz01}
  \begin{array}{c}
\Tr _{_{\!\! (n+1)}} \! \Bigl( \!
\tilde{T}_{k-1}(\frac{z_k}{z_{k-1}}) \cdots \tilde{T}_1(\frac{z_k}{z_1}) \, y_1(z_k)
\cdot \tilde{T}_{1}(z_k z_1) \cdots \tilde{T}_{k-1} (z_k z_{k-1}) \cdot
 \\ [0.3cm]
\cdot
\tilde{T}_{k}(z_k z_{k+1})   \cdots
   \tilde{T}_{n-1}(z_k z_{n}) \tilde{T}_n(z_k^2) \tilde{T}_{n-1}(\frac{z_k}{z_{n}})
 \cdots \tilde{T}_{k+1} (\frac{z_k}{z_{k+2}}) \tilde{T}_{k} (\frac{z_k}{z_{k+1}}) \! \Bigr)  \; .
\end{array}
 \ee
 Now we use relations (\ref{qtr1}) to obtain
 $$
 \Tr _{_{\!\! (n+1)}} \! \Bigl( \!  \tilde{T}_n(z^2) \! \Bigr) =
 \frac{(q^2 -z^2 \nu^2)(z^2 \nu +q^{-1})}{(z^2\nu +q)}
 \equiv N(z^2) \; .
 $$
Then for (\ref{qkz01}) we deduce
 \be
 \lb{qkz02}
 \begin{array}{c}
 B_k(\vec{z}) =  N(z_k^2) \,
 \Bigl( \! T_{k-1}(\frac{z_k}{z_{k-1}}) \cdots T_1(\frac{z_k}{z_1}) \,
 \cdot y_1(z_k)  \, T_{1}(z_k z_1) \cdots T_{k-1} (z_k z_{k-1}) \! \Bigr) \
  \cdot \\ [0.3cm]  \cdot  \Bigl( \! T_{k}(z_k z_{k+1})   \cdots
  T_{n-1}(z_k z_{n}) \cdot T_{n-1}(\frac{z_k}{z_{n}})
  \cdots    T_{k} (\frac{z_k}{z_{k+1}}) \! \Bigr)  = \\ [0.3cm] =
   N(z_k^2) \,  \tilde{y}_k (z_k, \vec{z}_{(1,k-1)})
  \cdot \overline{y}_k (z_k, \vec{z}_{(k+1,n)}) = 
  N(z_k^2) \, \Bigl( \prod\limits_{i=1}^{k-1}
\frac{(z_k q^2 -z_i)  \,  (z_i q^2  - z_k)}{q^2 z_k z_i} \Bigr)
 {\sf A}'_k(\vec{z}\,)  \; ,
\end{array}
 \ee
 where
   \be
   \lb{Bax05}
 {\sf A}'_k(\vec{z}\,) =
\tilde{y}_{k}'(z_k;\vec{z}_{(k-1)}) \cdot
 \overline{y}_k (z_k, \vec{z}_{(k+1,n)}) \; ,
 \ee
 $$
 \begin{array}{c}
 \overline{y}_k (z_k, \vec{z}_{(k+1,n)}) =
 \! \tilde{T}_{k}(z_k z_{k+1})   \cdots
   \tilde{T}_{n-1}(z_k z_{n}) \cdot \tilde{T}_{n-1}(\frac{z_k}{z_{n}})
 \cdots \tilde{T}_{k+1} (\frac{z_k}{z_{k+2}}) \tilde{T}_{k} (\frac{z_k}{z_{k+1}}) \; ,
 \end{array}
 $$
 and elements $\tilde{y}_k (x, \vec{z}_{(1,k-1)})$,
 $\tilde{y}_k'(x, \vec{z}_{(1,k-1)})$  were defined
 in (\ref{xxz55w}), (\ref{xxz55y}).

  Operators (\ref{qkz02}) are equal
  to the transfer-matrix operator
 $\tau_{n}(x ;\vec{z}_{(n)})$
 evaluated at the points $x=z_k$. Thus, by definition
 the operators $\{ B_1(\vec{z}\,) , \dots , B_n(\vec{z}\,) \}$
  form a commutative set of elements in the algebra $BMW_{n}(C)$:
 \be
 \lb{Bax06}
 [ B_k(\vec{z}\,) \, , \,  B_r(\vec{z}\,) ] = 0
 \;\;\;\;\; (\forall k,r = 1, \dots, n) \; .
 \ee
 Thus, operators $\{ B_1(\vec{z}\,) , \dots , B_n(\vec{z}\,) \}$
 for fixed parameters $\{ z_1, \dots, z_n \}$
 can be considered as generators of the Bethe subalgebra in $BMW_{n}(C)$.

 The validity of the identities (\ref{Bax06}) can be shown in different way.
 For this, in view of (\ref{qkz02}),  we need to prove the commutativity
 of the set of elements ${\sf A}'_k(\vec{z}\,) \; \in \; BMW_{n}(C)$
 which can be interpreted as analogs of flat connections (\ref{Zam77})
 for quantum Knizhnik-Zamolodchikov equations.
Taking into account (\ref{Zam10D}) we see that the map $\tilde{\rho}_{c}$: $B_n(C) \; \to \; BMW_{n}(C)$:
 \be
 \lb{Zam10K}
 \begin{array}{c}
 \tilde{\rho}_{c}(T_i) \; = \;  s_i \, \tilde{T}_{i}(z_i,z_{i+1})  \;\;\;\; (i=1,\dots,n-1) \; , \;\;\;
 \tilde{\rho}_{c}(T_0) \; = \;   y_1(z_1) \, s_0  \; ,
 \end{array}
 \ee
 where $s_0$ is defined in (\ref{Affbg1}) with $\sigma(z_1) = 1/z_1$ ,
 is the representation of $B_n(C)$. Then we have
 the following statement.

 \vspace{0.1cm}

 \noindent
 {\bf Proposition \ref{sec6}.4} {\it The flat connections
 ${\sf A}_i'(\vec{z}\,)$ (\ref{Bax05}) are images $\tilde{\rho}_{c}({\sf J}_i)$
 of the pairwise commuting elements
 $$
 {\sf J}_i = (T_{i-1}^{-1} \cdots T_1^{-1} T_0 T_1 \cdots T_{i-1})
 (T_i \cdots T_{n-1}  \cdot T_{n-1} \cdots T_i) \in B_n(C) \; ,
 $$
 which are obtained by the projection $T_n \to 1$ from the elements
 $J_i \in B_n(C^{(1)})$ given in (\ref{jucys1})}.  \\
 {\bf Proof.} The formula ${\sf A}_i'(\vec{z}\,) =
 \tilde{\rho}_{c}({\sf J}_i)$ can be checked directly with the use
 of definition (\ref{Bax05}) of ${\sf A}_i'(\vec{z}\,)$  and
 formulas (\ref{Zam10K}) for the map $\tilde{\rho}_{c}$. \hfill \qed

 \vspace{0.2cm}

\noindent
 {\bf Remark.} Using the
 special limit in (\ref{Zam77}), one can
 deduce the BMW analog of the Cherednik's connections
 \be
 \lb{qkz03}
A_k(\vec{z}) =
T_{k-1}\Bigl( \frac{ z_k}{z_{k-1}}\Bigr)
 \cdots T_1\Bigl(\frac{ z_k}{z_1}\Bigr)  \cdot
 y_1^{-1} {\sf T}_1^{-1}  \cdots {\sf T}_{n-1}^{-1} \, D_{z_k}  \cdot
 T_{n-1}\Bigl(\frac{z_k}{z_{n}}\Bigr)
 \cdots
  T_{k} \Bigl( \frac{z_k}{z_{k+1}}
  \! \Bigr) \; \in \; BMW_n(C) \; ,
 \ee
 which were presented for the Hecke algebra case in \cite{Ch} (see there
 eq. (4.12) in Section 4.2).
 The finite difference operator $D_{z_k}$  is given in (\ref{Zam57})
 with $\tilde{\bar{z}}= \frac{c}{c'} z$.
 In (\ref{qkz03}) we have to take into account that Cherednik's
 affine elements $Y_k$ are related to ours
 by $Y_k = y_k^{-1}$.

 To rewrite our expression (\ref{Zam77}) to the Cherednik's one
  (\ref{qkz03}) we need to convert the factor
  \be
  \lb{qkz05}
  L_1(z_k)  \, T_{1}( c z_k z_1) \cdots T_{k-1} (c z_k z_{k-1})
  \cdot   T_{k}(c z_k z_{k+1})   \cdots T_{n-1}(c z_k z_{n})
  \, \bar{L}_n( \frac{1}{c z_k} ) \; ,
  \ee
  entered into the expression (\ref{Zam77}) to
  the factor $y_1^{-1} {\sf T}_1^{-1}  \cdots {\sf T}_{n-1}^{-1}$.
  It can be done if
  we first make in (\ref{Zam77}) the redefinition
  of all spectral parameters $z_r \to t z_r$ and then
  consider the limit $t \to \infty$. To do this we note
  that only the product (\ref{qkz05}) in (\ref{Zam77})
  will be dependent on $t$, where we have to use limits
  $$
  \lim_{t \to \infty} T_r(t^2 c z_k z_r) =
  {\sf T}_r - (q - q^{-1})  + (q - q^{-1})  \kappa_r = {\sf T}_r^{-1} \; ,
  $$
  $$
  \lim_{t \to \infty} L_1(t z_k)    = \frac{1}{c} \, y_1^{-1}
  = \frac{1}{c} \, {\sf T}_0^{-1}\; ,
\;\;\;\;
 \lim_{t \to \infty} \bar{L}_n \Bigl( \frac{1}{t c z_k} \Bigr)
 =  \bar{y}_n = {\sf T}_n \; .
  $$
  Here we used the expressions for baxterized elements
  (\ref{a00}), (\ref{rea15}) and (\ref{rea151}).
  Thus, the limit of the factor (\ref{qkz05}) is
  $$
  y_1^{-1} \cdot {\sf T}_1^{-1}  \cdots {\sf T}_{n-1}^{-1} \cdot \bar{y}_n =
   {\sf T}_0^{-1} \cdot {\sf T}_1^{-1}  \cdots {\sf T}_{n-1}^{-1}
   \cdot {\sf T}_n \equiv {\sf X} \; ,
   $$
   and for the limit of the connection (\ref{Zam77})
   we obtain expression
 $$
 A_k'(\vec{z}) =
T_{k-1}\Bigl( \frac{ z_k}{z_{k-1}}\Bigr)
 \cdots T_1\Bigl(\frac{ z_k}{z_1}\Bigr)  \cdot
 {\sf X} \,  D_{z_k} \cdot
 T_{n-1}\Bigl(\frac{z_k}{z_{n}}\Bigr)
 \cdots T_{k+1} \Bigl(\frac{z_k}{z_{k+2}}\Bigr)  T_{k} \Bigl( \frac{z_k}{z_{k+1}}
  \! \Bigr)  \; \in \; BMW_n(C^{(1)})\; ,
 $$
 which is a generalization of (\ref{qkz03}).
 The projection ${\sf T}_n \to 1$ for connection $A_k'(\vec{z})$ gives
 the BMW analog of the Cherednik's connection (\ref{qkz03}).

\vspace{1cm}

{\bf Acknowledgment.} The work of API was supported by RSCF grant
14-11-00598.

\end{document}